# A recurrence relation for the Li/Keiper constants in terms of the Stieltjes constants

Donal F. Connon

dconnon@btopenworld.com

10 February 2009

**Abstract**

A recurrence relation for the Li/Keiper constants in terms of the Stieltjes constants is derived in this paper.

In addition, we also report a formula for the Stieltjes constants in terms of the higher derivatives, $\varsigma^{(n)}(0)$, of the Riemann zeta function evaluated at $s=0$. A formula for the Stieltjes constants in terms of the (exponential) complete Bell polynomials containing the eta constants $\eta_n$ as the arguments is also derived.

## 1. Introduction

The Riemann xi function $\xi(s)$ is defined as

(1.1) $\qquad \xi(s) = \frac{1}{2} s(s-1) \pi^{-s/2} \Gamma(s/2) \varsigma(s)$

and we see from the functional equation for the Riemann zeta function $\varsigma(s)$

(1.2) $\qquad \varsigma(1-s) = 2(2\pi)^{-s} \Gamma(s) \cos(\pi s/2) \varsigma(s)$

that $\xi(s)$ satisfies the functional equation

(1.3) $\qquad \xi(s) = \xi(1-s)$

In 1996, Li [24] defined the sequence of numbers $(\lambda_n)$ by

(1.4) $\qquad \lambda_n = \dfrac{1}{(n-1)!} \dfrac{d^n}{ds^n} [s^{n-1} \log \xi(s)] \bigg|_{s=1}$

and proved that a necessary and sufficient condition for the non-trivial zeros $\rho$ of the Riemann zeta function to lie on the critical line $s = \dfrac{1}{2} + i\tau$ is that $\lambda_n$ is non-negative for every positive integer $n$. Earlier in 1991, Keiper [21] showed that if the Riemann hypothesis is true, then $\lambda_n > 0$ for all $n \geq 1$.

Li also showed that

(1.5) $\qquad \lambda_n = \sum_{\rho} \left[ 1 - \left(1 - \dfrac{1}{\rho}\right)^n \right]$

Taking logarithms of (1.1) and noting that $\dfrac{s}{2} \Gamma\left(\dfrac{s}{2}\right) = \Gamma\left(1 + \dfrac{s}{2}\right)$ we see that

(1.6) $\qquad \log \xi(s) = \log \Gamma\left(1 + \dfrac{s}{2}\right) - \dfrac{s}{2} \log \pi + \log[(s-1)\varsigma(s)]$

We also have the Maclaurin expansion about $s = 1$

(1.7) $\qquad \log \xi(s) = -\log 2 - \sum_{k=1}^{\infty} (-1)^k \dfrac{\sigma_k}{k} (s-1)^k$

where the constant term in (1.7) arises because

$$\lim_{s \to 1} \xi(s) = \dfrac{1}{2} \pi^{-\frac{1}{2}} \Gamma(1/2) \lim_{s \to 1} [(s-1)\varsigma(s)] = \dfrac{1}{2}$$

We note that the coefficients $\sigma_k$ are defined by

$$\sigma_k = \frac{(-1)^{k+1}}{(k-1)!} \frac{d^k}{ds^k} \log \xi(s) \bigg|_{s=1}$$

and we see that

$$\frac{d}{ds} \log \xi(s) = \frac{\xi'(s)}{\xi(s)} = -\sum_{k=1}^{\infty} (-1)^k \sigma_k (s-1)^{k-1}$$

and comparing this with (1.4) we immediately see that $\lambda_1 = \sigma_1$.

Differentiating (1.6) results in

(1.8) $$\frac{d}{ds} \log \xi(s) = \frac{1}{2}\psi\left(1+\frac{s}{2}\right) - \frac{1}{2}\log \pi + \frac{\frac{d}{ds}[(s-1)\varsigma(s)]}{(s-1)\varsigma(s)}$$

which, referring to (1.4), is clearly relevant in computing $\lambda_1$.

Before proceeding any further we need to recall some details regarding the Stieltjes constants.

## 2. A brief survey of the Stieltjes constants

The Stieltjes constants $\gamma_n(u)$ are the coefficients in the Laurent expansion of the Hurwitz zeta function $\varsigma(s,u)$ about $s=1$

(2.1) $$\varsigma(s,u) = \sum_{n=0}^{\infty} \frac{1}{(n+u)^s} = \frac{1}{s-1} + \sum_{n=0}^{\infty} \frac{(-1)^n}{n!} \gamma_n(u)(s-1)^n$$

and $\gamma_0(u) = -\psi(u)$, where $\psi(u)$ is the digamma function which is the logarithmic derivative of the gamma function $\psi(u) = \frac{d}{du} \log \Gamma(u)$. It is easily seen from the definition of the Hurwitz zeta function that $\varsigma(s,1) = \varsigma(s)$ and accordingly that $\gamma_n(1) = \gamma_n$. Further details of the Stieltjes constants are contained in [5] to [13] inclusive.

The generalised Euler-Mascheroni constants $\gamma_n$ (or Stieltjes constants) are the coefficients of the Laurent expansion of the Riemann zeta function $\varsigma(s)$ about $s=1$

(2.2) $$\varsigma(s) = \frac{1}{s-1} + \sum_{n=0}^{\infty} \frac{(-1)^n}{n!} \gamma_n (s-1)^n$$

Since $\lim_{s \to 1}\left[\varsigma(s) - \frac{1}{s-1}\right] = \gamma$ it is clear that $\gamma_0 = \gamma$. It may be shown, as in [20, p.4], that

$$(2.3) \quad \gamma_n = \lim_{N \to \infty}\left[\sum_{k=1}^{N} \frac{\log^n k}{k} - \frac{\log^{n+1} N}{n+1}\right] = \lim_{N \to \infty}\left[\sum_{k=1}^{N} \frac{\log^n k}{k} - \int_1^N \frac{\log^n t}{t} dt\right]$$

where, throughout this paper, we define $\log^0 1 = 1$.

We now refer to the Hasse identity [19] for the Hurwitz zeta function $\varsigma(s,u)$ which is valid for all $s \in \mathbb{C}$ provided $\text{Re}(s) \neq 1$

$$(2.4) \quad \varsigma(s,u) = \frac{1}{s-1} \sum_{i=0}^{\infty} \frac{1}{i+1} \sum_{j=0}^{i} \binom{i}{j} \frac{(-1)^j}{(u+j)^{s-1}}$$

and we have the well-known limit

$$(2.5) \quad \lim_{s \to 1}(s-1)\varsigma(s,u) = \sum_{i=0}^{\infty} \frac{1}{i+1} \sum_{j=0}^{i} \binom{i}{j}(-1)^j = \sum_{i=0}^{\infty} \frac{1}{i+1} \delta_{i,0} = 1$$

where $\delta_{i,j}$ is the Kronecker delta symbol.

We see from (2.4) that

$$(2.6) \quad \frac{d^{n+1}}{ds^{n+1}}[(s-1)\varsigma(s,u)] = (-1)^{n+1} \sum_{i=0}^{\infty} \frac{1}{i+1} \sum_{j=0}^{i} \binom{i}{j} \frac{(-1)^j \log^{n+1}(u+j)}{(u+j)^{s-1}}$$

and thus

$$(2.7) \quad \frac{d^{n+1}}{ds^{n+1}}[(s-1)\varsigma(s,u)]\bigg|_{s=1} = (-1)^{n+1} \sum_{i=0}^{\infty} \frac{1}{i+1} \sum_{j=0}^{i} \binom{i}{j}(-1)^j \log^{n+1}(u+j)$$

We previously showed in [15] that

$$(2.8) \quad \gamma_n(u) = -\frac{1}{n+1} \sum_{i=0}^{\infty} \frac{1}{i+1} \sum_{j=0}^{i} \binom{i}{j}(-1)^j \log^{n+1}(u+j)$$

and comparing this with (2.7) allows us to conclude for $n \geq 0$ that

$$(2.9) \quad \frac{d^{n+1}}{ds^{n+1}}[(s-1)\varsigma(s,u)]\bigg|_{s=1} = (-1)^n (n+1)\gamma_n(u)$$

This may also be more directly obtained by differentiating (2.1) $p+1$ times where we see that

$$\frac{d^{p+1}}{ds^{p+1}}[(s-1)\varsigma(s,u)] = \sum_{n=0}^{\infty}\frac{(-1)^n}{n!}\gamma_n(u)(n+1)n(n+1)\cdots(n+1-p)(s-1)^{n-p}$$

With $u=1$ in (2.9) we have

(2.10) $$\left.\frac{d^{n+1}}{ds^{n+1}}[(s-1)\varsigma(s)]\right|_{s=1} = (-1)^n(n+1)\gamma_n$$

We now evaluate (1.8) at $s=1$ and obtain

$$\lambda_1 = \left.\frac{d}{ds}[\log\xi(s)]\right|_{s=1} = \frac{1}{2}\psi\left(\frac{3}{2}\right) - \frac{1}{2}\log\pi + \left.\frac{d}{ds}[(s-1)\varsigma(s)]\right|_{s=1}$$

Using [30, p.20]

$$\psi\left(n+\frac{1}{2}\right) = -\gamma - 2\log 2 + 2\sum_{k=0}^{n-1}\frac{1}{2k+1}$$

we see that

(2.11) $$\psi\left(\frac{3}{2}\right) = -\gamma - 2\log 2 + 2$$

From (2.10) we have

(2.12) $$\left.\frac{d}{ds}[(s-1)\varsigma(s)]\right|_{s=1} = \gamma_0 = \gamma$$

and we therefore easily compute the first Li/Keiper constant

(2.13) $$\lambda_1 = \left.\frac{d}{ds}[\log\xi(s)]\right|_{s=1} = \frac{1}{2}\psi\left(\frac{3}{2}\right) + \gamma - \frac{1}{2}\log\pi = -\frac{1}{2}\log\pi + \frac{1}{2}\gamma + 1 - \log 2$$

We note from (2.11) that $\psi\left(\frac{3}{2}\right)$ is positive; in addition we see that $\gamma - \frac{1}{2}\log\pi > 0$ and therefore $\lambda_1$ is positive. Its approximate value is [5]

(2.14) $$\sigma_1 = \lambda_1 \simeq 0.023...$$

However, in order to determine $\lambda_n$ we would need to calculate $\dfrac{d^n}{ds^n}[s^{n-1}\log \xi(s)]$ and reference to (1.6) shows that this in turn would require us to compute

$$\frac{d^n}{ds^n}[s^{n-1}\log[(s-1)\varsigma(s)]] = \frac{d^{n-1}}{ds^{n-1}}\left[(n-1)s^{n-2}\log[(s-1)\varsigma(s)] + s^{n-1}\frac{\dfrac{d}{ds}[(s-1)\varsigma(s)]}{(s-1)\varsigma(s)}\right]$$

It is remarkable that determining the veracity, or otherwise, of the most famous hypothesis in mathematics is so dependent upon obtaining a "manageable" formula for the higher derivatives of the quotient of two functions.

Formulae for the Li/Keiper constants have been given by Bombieri and Lagarias [4], Maślanka [26] and Coffey [13], the latter being based on an identity determined by Matsuoka [27]. However, it is not easy to discern particular values of $\lambda_n$ from these representations.

It is because of this inherent difficulty that we now proceed to go for second best and derive a recurrence relationship for the $\lambda_n$ constants.

**3. A recurrence relation for the Li/Keiper constants**

We recall (1.7) above

$$\log \xi(s) = -\log 2 - \sum_{k=1}^{\infty}(-1)^k \frac{\sigma_k}{k}(s-1)^k$$

and since $\log \xi(s) = \log \xi(1-s)$ we see that

$$\log \xi(s) = -\log 2 - \sum_{k=1}^{\infty}\frac{\sigma_k}{k}s^k$$

We therefore have from (1.6)

(3.1) $\quad \log \Gamma\left(1+\dfrac{s}{2}\right) + \log 2 - \dfrac{s}{2}\log \pi + \log[(s-1)\varsigma(s)] = -\sum_{k=1}^{\infty}\dfrac{\sigma_k}{k}s^k$

and differentiation gives us

(3.2) $\quad \dfrac{1}{2}\psi\left(1+\dfrac{s}{2}\right) - \dfrac{1}{2}\log \pi + \dfrac{\dfrac{d}{ds}[(s-1)\varsigma(s)]}{(s-1)\varsigma(s)} = -\sum_{k=1}^{\infty}\dfrac{\sigma_k}{k}ks^{k-1}$

We now multiply (3.2) across by $(s-1)\varsigma(s)$ to obtain

(3.3)

$$\frac{1}{2}(s-1)\varsigma(s)\psi\left(1+\frac{s}{2}\right) - \frac{1}{2}\log \pi (s-1)\varsigma(s) + \frac{d}{ds}[(s-1)\varsigma(s)] = -\sum_{k=1}^{\infty}\frac{\sigma_k}{k}ks^{k-1}[(s-1)\varsigma(s)]$$

and with $s=1$ we have

$$\frac{1}{2}\psi\left(\frac{3}{2}\right) - \frac{1}{2}\log \pi + \gamma = -\sum_{k=1}^{\infty}\frac{\sigma_k}{k}k$$

Using

$$\log \xi(s) = -\log 2 - \sum_{k=1}^{\infty}\frac{\sigma_k}{k}s^k$$

we see that

$$\frac{d^n}{ds^n}[s^{n-1}\log \xi(s)] = -\sum_{k=1}^{\infty}\frac{\sigma_k}{k}(k+n-1)\cdots k \, s^{k+n-1}$$

and hence referring to (1.4) we have

$$\lambda_n = \frac{1}{(n-1)!}\frac{d^n}{ds^n}[s^{n-1}\log \xi(s)]\bigg|_{s=1} = -\frac{1}{(n-1)!}\sum_{k=1}^{\infty}\frac{\sigma_k}{k}(k+n-1)\cdots k$$

We then see that for $n \geq 1$

(3.4)    $$(n-1)!\lambda_n = -\sum_{k=1}^{\infty}\frac{\sigma_k}{k}k\cdots(k+n-1)$$

and for example we have

$$\lambda_1 = -\frac{1}{0!}\sum_{k=1}^{\infty}\frac{\sigma_k}{k}k = -\sum_{k=1}^{\infty}\sigma_k = \sigma_1$$

$$\lambda_2 = -\frac{1}{1!}\sum_{k=1}^{\infty}\frac{\sigma_k}{k}k(k+1)$$

$$\lambda_3 = -\frac{1}{2!}\sum_{k=1}^{\infty}\frac{\sigma_k}{k}k(k+1)(k+2)$$

We then see that

$$\frac{1}{2}\psi\left(\frac{3}{2}\right)+\gamma-\frac{1}{2}\log\pi = \lambda_1$$

which is the same as (2.11).

Differentiating equation (3.3) gives us

$$(s-1)\varsigma(s)\frac{1}{2^2}\psi^{(1)}\left(1+\frac{s}{2}\right)+\frac{1}{2}\psi\left(1+\frac{s}{2}\right)\frac{d}{ds}[(s-1)\varsigma(s)]-\frac{1}{2}\log\pi\frac{d}{ds}[(s-1)\varsigma(s)]+\frac{d^2}{ds^2}[(s-1)\varsigma(s)]$$

$$=-\sum_{k=1}^{\infty}\frac{\sigma_k}{k}ks^{k-1}\frac{d}{ds}[(s-1)\varsigma(s)]-\sum_{k=1}^{\infty}\frac{\sigma_k}{k}k(k-1)s^{k-2}[(s-1)\varsigma(s)]$$

and with $s=1$ we have

$$\frac{1}{4}\psi^{(1)}\left(\frac{3}{2}\right)+\frac{1}{2}\psi\left(\frac{3}{2}\right)\gamma-\frac{1}{2}\gamma\log\pi-2\gamma_1 = -\gamma\sum_{k=1}^{\infty}\sigma_k-\sum_{k=1}^{\infty}\frac{\sigma_k}{k}k(k-1)$$

We have $k(k-1) = k(k+1)-2k$ and hence

$$-\gamma\sum_{k=1}^{\infty}\sigma_k-\sum_{k=1}^{\infty}\frac{\sigma_k}{k}k(k-1) = -\gamma\sum_{k=1}^{\infty}\frac{\sigma_k}{k}k-\sum_{k=1}^{\infty}\frac{\sigma_k}{k}k(k+1)+2\sum_{k=1}^{\infty}\frac{\sigma_k}{k}k$$

and, therefore, using (3.4), we have

(3.5) $$\frac{1}{4}\psi^{(1)}\left(\frac{3}{2}\right)+\frac{1}{2}\gamma\psi\left(\frac{3}{2}\right)-\frac{1}{2}\gamma\log\pi-2\gamma_1 = (\gamma-2)\lambda_1+\lambda_2$$

Since $\psi(1+x) = \psi(x)+\frac{1}{x}$ we see that

$$\psi^{(n)}(1+x) = \psi^{(n)}(x)+(-1)^n n! x^{-n-1}$$

and with $x = 1/2$ we have

$$\psi^{(n)}(3/2) = \psi^{(n)}(1/2)+(-1)^n n! 2^{n+1}$$

We have the well-known result from [30]

$$\psi^{(n)}(1/2) = (-1)^{n+1} n! [2^{n+1}-1]\varsigma(n+1)$$

and hence we have

(3.5.1) $$\psi^{(n)}(3/2) = (-1)^{n+1} n! \left( \left[ 2^{n+1} - 1 \right] \varsigma(n+1) - 2^{n+1} \right)$$

In fact we have from [30, p.22]

(3.5.2) $$\psi^{(n)}(x) = (-1)^{n+1} n! \varsigma(n+1, x)$$

and hence we deduce that for $n \geq 1$

(3.5.3) $$\psi^{(n)}(3/2) = (-1)^{n+1} c_n \quad \text{where } c_n > 0$$

We thus obtain

$$\psi^{(1)}\left(\frac{3}{2}\right) = 3\varsigma(2) - 4$$

We have from (2.11)

$$\psi\left(\frac{3}{2}\right) = -\gamma - 2\log 2 + 2$$

and from (2.13)

$$\lambda_1 = -\frac{1}{2}\log \pi + \frac{1}{2}\gamma + 1 - \log 2$$

Therefore we have the second Li/Keiper constant

(3.6) $$\lambda_2 = \frac{3}{4}\varsigma(2) + 1 + \gamma - \gamma^2 - 2\log 2 - \log \pi - 2\gamma_1$$

Differentiating equation (3.3) $n+1$ times using the Leibniz rule gives us

(3.7)
$$\sum_{m=0}^{n+1} \binom{n+1}{m} \frac{d^m}{ds^m}[(s-1)\varsigma(s)] \frac{1}{2^{n+2-m}} \psi^{(n+1-m)}\left(1 + \frac{s}{2}\right) - \frac{1}{2}\log \pi \frac{d^{n+1}}{ds^{n+1}}[(s-1)\varsigma(s)]$$

$$+ \frac{d^{n+2}}{ds^{n+2}}[(s-1)\varsigma(s)] = -\sum_{k=1}^{\infty} \frac{\sigma_k}{k} \sum_{m=0}^{n+1} \binom{n+1}{m} \frac{d^m}{ds^m}[ks^{k-1}] \frac{d^{n+1-m}}{ds^{n+1-m}}[(s-1)\varsigma(s)]$$

First of all, isolating the term for $m = 0$ and evaluating at $s = 1$ we have since $\lim_{s \to 1}[(s-1)\varsigma(s)] = 1$

$$S_1 = \sum_{m=0}^{n+1} \binom{n+1}{m} \frac{d^m}{ds^m}[(s-1)\varsigma(s)] \frac{1}{2^{n+2-m}} \psi^{(n+1-m)}\left(1 + \frac{s}{2}\right)\bigg|_{s=1}$$

$$= \frac{1}{2^{n+2}} \psi^{(n+1)}\left(1 + \frac{s}{2}\right) + \sum_{m=1}^{n+1} \binom{n+1}{m} \frac{d^m}{ds^m}[(s-1)\varsigma(s)] \frac{1}{2^{n+2-m}} \psi^{(n+1-m)}\left(1 + \frac{s}{2}\right)\bigg|_{s=1}$$

and then rebasing the index to $m = p + 1$ results in

$$S_1 = \frac{1}{2^{n+2}} \psi^{(n+1)}\left(1 + \frac{s}{2}\right) + \sum_{p=0}^{n} \binom{n+1}{p+1} \frac{d^{p+1}}{ds^{p+1}}[(s-1)\varsigma(s)] \frac{1}{2^{n+1-p}} \psi^{(n-p)}\left(1 + \frac{s}{2}\right)\bigg|_{s=1}$$

Noting from (2.9) that

$$\frac{d^{p+1}}{ds^{p+1}}[(s-1)\varsigma(s)]\bigg|_{s=1} = (-1)^p (p+1)\gamma_p$$

we obtain

$$S_1 = \frac{1}{2^{n+2}} \psi^{(n+1)}\left(\frac{3}{2}\right) + \sum_{p=0}^{n} \binom{n+1}{p+1} (-1)^p (p+1)\gamma_p \frac{1}{2^{n+1-p}} \psi^{(n-p)}\left(\frac{3}{2}\right)$$

Carrying out similar operations for the summation on the right-hand side of (3.7), where this time we isolate the term for $m = n + 1$, we get

$$S_2 = -\sum_{k=1}^{\infty} \frac{\sigma_k}{k} \sum_{m=0}^{n+1} \binom{n+1}{m} \frac{d^m}{ds^m}[ks^{k-1}] \frac{d^{n+1-m}}{ds^{n+1-m}}[(s-1)\varsigma(s)]\bigg|_{s=1}$$

$$= -\sum_{k=1}^{\infty} \frac{\sigma_k}{k} \frac{d^{n+1}}{ds^{n+1}}[ks^{k-1}] - \sum_{k=1}^{\infty} \frac{\sigma_k}{k} \sum_{m=0}^{n} \binom{n+1}{m} \frac{d^m}{ds^m}[ks^{k-1}] \frac{d^{n+1-m}}{ds^{n+1-m}}[(s-1)\varsigma(s)]\bigg|_{s=1}$$

$$= -\sum_{k=1}^{\infty} \frac{\sigma_k}{k} k(k-1)\ldots(k-n-1) - \sum_{k=1}^{\infty} \frac{\sigma_k}{k} \sum_{m=0}^{n} \binom{n+1}{m} k(k-1)\ldots(k-m) \frac{d^{n+1-m}}{ds^{n+1-m}}[(s-1)\varsigma(s)]\bigg|_{s=1}$$

$$= -\sum_{k=1}^{\infty} \frac{\sigma_k}{k} k(k-1)\ldots(k-n-1) - \sum_{k=1}^{\infty} \frac{\sigma_k}{k} \sum_{m=0}^{n} \binom{n+1}{m} k(k-1)\ldots(k-m)(-1)^{n-m}(n-m+1)\gamma_{n-m}$$

Therefore, the evaluation of (3.7) at $s=1$ results in

$$\frac{1}{2^{n+2}}\psi^{(n+1)}\left(\frac{3}{2}\right)+\sum_{m=0}^{n}\binom{n+1}{m+1}(-1)^m(m+1)\gamma_m\frac{1}{2^{n+1-m}}\psi^{(n-m)}\left(\frac{3}{2}\right)$$

$$-\frac{1}{2}(-1)^n(n+1)\gamma_n\log\pi+(-1)^{n+1}(n+2)\gamma_{n+1}$$

$$=-\sum_{k=1}^{\infty}\frac{\sigma_k}{k}k(k-1)...(k-n-1)-\sum_{k=1}^{\infty}\frac{\sigma_k}{k}\sum_{m=0}^{n}\binom{n+1}{m}k(k-1)...(k-m)(-1)^{n-m}(n-m+1)\gamma_{n-m}$$

which may be written as

(3.8) $$\frac{1}{2^{n+2}}\psi^{(n+1)}\left(\frac{3}{2}\right)+\sum_{m=0}^{n}\binom{n+1}{m+1}(-1)^m(m+1)\gamma_m\frac{1}{2^{n+1-m}}\psi^{(n-m)}\left(\frac{3}{2}\right)$$

$$+\frac{1}{2}(-1)^{n+1}(n+1)\gamma_n\log\pi+(-1)^{n+1}(n+2)\gamma_{n+1}$$

$$=-\sum_{k=1}^{\infty}\frac{\sigma_k}{k}k(k-1)...(k-n-1)-\sum_{m=0}^{n}\binom{n+1}{m}(-1)^{n-m}(n-m+1)\gamma_{n-m}\sum_{k=1}^{\infty}\frac{\sigma_k}{k}k(k-1)...(k-m)$$

The right-hand side of this equation clearly requires a modicum of simplification which we now endeavour to provide.

The rising factorial $[x]^p$ is defined by

$$[x]^p = x(x+1)...(x+p-1)$$

and it is shown in [2, p.36] that

$$[x]^p = \sum_{j=1}^{p}\frac{p!}{j!}\binom{p-1}{j-1}[x]_j$$

where the falling factorial $[x]_j$ is defined by

$$[x]_j = x(x-1)...(x-j+1) = \sum_{i=0}^{j}s(j,i)x^i$$

and $s(j,i)$ are known as the Stirling numbers of the first kind. We therefore have

$$x(x+1)\ldots(x+p-1) = \sum_{j=1}^{p} \frac{p!}{j!}\binom{p-1}{j-1} x(x-1)\ldots(x-j+1)$$

Letting $x \to -x$ we see that

$$x(x-1)\ldots(x-p+1) = (-1)^p \sum_{j=1}^{p} \frac{p!}{j!}\binom{p-1}{j-1}(-1)^j x(x+1)\ldots(x+j-1)$$

and with $x = k$ we have

(3.9) $\quad k(k-1)\ldots(k-p+1) = (-1)^p \sum_{j=1}^{p} \frac{p!}{j!}\binom{p-1}{j-1}(-1)^j k(k+1)\ldots(k+j-1)$

We then have for the first term in the right-hand side of (3.8)

$$S_3 = -\sum_{k=1}^{\infty} \frac{\sigma_k}{k} k(k-1)\ldots(k-n-1) = -(-1)^n \sum_{j=1}^{n+2} \frac{(n+2)!}{j!}\binom{n+1}{j-1}(-1)^j \sum_{k=1}^{\infty} \frac{\sigma_k}{k} k(k+1)\ldots(k+j-1)$$

Using (3.4) for $n \geq 1$

$$(n-1)!\lambda_n = -\sum_{k=1}^{\infty} \frac{\sigma_k}{k} k \cdots (k+n-1)$$

we may then express $S_3$ in terms of the Li/Keiper constants as

(3.10) $\quad S_3 = (-1)^n \sum_{j=1}^{n+2} \frac{(n+2)!}{j!}\binom{n+1}{j-1}(-1)^j (j-1)!\lambda_j = (-1)^n \sum_{j=1}^{n+2} \frac{(n+2)!}{j}\binom{n+1}{j-1}(-1)^j \lambda_j$

We now consider the second term in the right-hand side of (3.8)

$$S_4 = -\sum_{m=0}^{n}\binom{n+1}{m}(-1)^{n-m}(n-m+1)\gamma_{n-m} \sum_{k=1}^{\infty} \frac{\sigma_k}{k} k(k-1)\ldots(k-m)$$

which on using (3.9) becomes

$$= -\sum_{m=0}^{n}\binom{n+1}{m}(-1)^{n-m}(n-m+1)\gamma_{n-m}(-1)^{m+1}\sum_{j=1}^{m+1}\frac{(m+1)!}{j!}\binom{m}{j-1}(-1)^j \sum_{k=1}^{\infty}\frac{\sigma_k}{k} k(k+1)\ldots(k+j-1)$$

$$= \sum_{m=0}^{n}\binom{n+1}{m}(-1)^{n-m}(n-m+1)\gamma_{n-m}(-1)^{m+1}\sum_{j=1}^{m+1}\frac{(m+1)!}{j!}\binom{m}{j-1}(-1)^j (j-1)!\lambda_j$$

(3.11)
$$= \sum_{m=0}^{n} \binom{n+1}{m}(-1)^{n-m}(n-m+1)\gamma_{n-m}(-1)^{m+1} \sum_{j=1}^{m+1} \frac{(m+1)!}{j}\binom{m}{j-1}(-1)^j \lambda_j$$

Using (3.10) and (3.11) we may then write (3.8) as

(3.12) $$\frac{1}{2^{n+2}}\psi^{(n+1)}\left(\frac{3}{2}\right) + \sum_{m=0}^{n}\binom{n+1}{m+1}(-1)^m(m+1)\gamma_m \frac{1}{2^{n+1-m}}\psi^{(n-m)}\left(\frac{3}{2}\right)$$

$$+ \frac{1}{2}(-1)^{n+1}(n+1)\gamma_n \log\pi + (-1)^{n+1}(n+2)\gamma_{n+1}$$

$$= (-1)^n \sum_{j=1}^{n+2} \frac{(n+2)!}{j}\binom{n+1}{j-1}(-1)^j \lambda_j$$

$$+ \sum_{m=0}^{n}\binom{n+1}{m}(-1)^{n-m}(n-m+1)\gamma_{n-m}(-1)^{m+1}\sum_{j=1}^{m+1}\frac{(m+1)!}{j}\binom{m}{j-1}(-1)^j \lambda_j$$

After some algebra and using the following elementary binomial identities

$$\binom{n+1}{m+1}(m+1) = (n+1)\binom{n}{m} \qquad \binom{n+1}{m}(n-m+1) = (m+1)\binom{n}{m+1}$$

$$\frac{1}{j}\binom{n+1}{j-1} = \frac{1}{n+2}\binom{n+2}{j} \qquad \frac{1}{j}\binom{m}{j-1} = \frac{1}{m+1}\binom{m+1}{j}$$

we have the required recurrence relation

(3.13) $$\frac{1}{2^{n+2}}\psi^{(n+1)}\left(\frac{3}{2}\right) + \frac{(n+1)}{2^{n+1}}\sum_{m=0}^{n}\binom{n}{m}(-1)^m \gamma_m 2^m \psi^{(n-m)}\left(\frac{3}{2}\right)$$

$$+ \frac{1}{2}(-1)^{n+1}(n+1)\gamma_n \log\pi + (-1)^{n+1}(n+2)\gamma_{n+1}$$

$$= (-1)^n(n+1)!\sum_{j=1}^{n+2}\binom{n+2}{j}(-1)^j \lambda_j + (-1)^{n+1}\sum_{m=0}^{n}\binom{n}{m+1}(m+1)!\gamma_{n-m}\sum_{j=1}^{m+1}\binom{m+1}{j}(-1)^j \lambda_j$$

For example with $n=0$ we get

(3.14) $$\frac{1}{4}\psi^{(1)}\left(\frac{3}{2}\right)+\frac{1}{2}\gamma\psi\left(\frac{3}{2}\right)-\frac{1}{2}\gamma\log\pi-2\gamma_1 = -2\lambda_1+\lambda_2+\gamma\lambda_1$$

which was previously derived in (3.5). Using (2.13) this may be expressed as

$$\frac{1}{4}\psi^{(1)}\left(\frac{3}{2}\right)-\gamma^2-2\gamma_1+2\lambda_1 = \lambda_2$$

and referring to (4.4) this becomes

$$\frac{1}{4}\psi^{(1)}\left(\frac{3}{2}\right)-\eta_1+2\lambda_1 = \lambda_2$$

and we note that $\psi^{(1)}\left(\frac{3}{2}\right)$ and $\eta_1$ are both positive (however at this stage we still need to substitute numerical values of the various constants to determine that $\lambda_2$ is non-negative).

□

In what follows we show how it is possible to determine the signs of the first few $\lambda_n$ with a modicum of numerical calculations.

Using the Jacobi theta function, Coffey [5] showed in 2003 that even derivatives of $\xi(s)$ are positive for all real values of $s$, while odd derivatives are positive for $s \geq \frac{1}{2}$ and negative for $s < \frac{1}{2}$. In particular, $\xi^{(n)}(1)$ is positive for $n \geq 1$. This was also proved by Freitas [18] in a different manner in 2005.

We note from (1.7) that

$$\xi'(s) = -\xi(s)\sum_{k=1}^{\infty}(-1)^k\frac{\sigma_k}{k}k(s-1)^{k-1} = -\xi(s)\sum_{k=1}^{\infty}\frac{\sigma_k}{k}ks^{k-1}$$

Since $\xi(1) = \xi(0) = \frac{1}{2}$ this results in

(3.15) $$\xi'(1) = \frac{1}{2}\sigma_1 = \frac{1}{2}\lambda_1$$

and we immediately deduce that both $\lambda_1$ and $\sigma_1$ are positive.

The (odd) derivative $\xi'(0) = \frac{1}{2}\sum_{k=1}^{\infty}\sigma_k = -\frac{1}{2}\lambda_1$ which again shows that $\lambda_1$ is positive and that $\sum_{k=1}^{\infty}\sigma_k$ is negative.

We have the second derivative

$$\xi^{(2)}(s) = -\xi(s)\sum_{k=1}^{\infty}\frac{\sigma_k}{k}k(k-1)s^{k-2} - \xi^{(1)}(s)\sum_{k=1}^{\infty}\frac{\sigma_k}{k}ks^{k-1}$$

and thus

$$\xi^{(2)}(1) = -\xi(1)\sum_{k=1}^{\infty}\frac{\sigma_k}{k}k(k-1) - \xi^{(1)}(1)\sum_{k=1}^{\infty}\frac{\sigma_k}{k}k$$

Since $k(k-1) = k(k+1) - 2k$ we have

$$\xi^{(2)}(1) = -\xi(1)\sum_{k=1}^{\infty}\frac{\sigma_k}{k}k(k+1) + 2\xi(1)\sum_{k=1}^{\infty}\frac{\sigma_k}{k}k - \xi^{(1)}(1)\sum_{k=1}^{\infty}\frac{\sigma_k}{k}k$$

and hence, using (3.4) and (3.15), we obtain

(3.16) $\quad \xi^{(2)}(1) = \frac{1}{2}\lambda_2 - \lambda_1 + \frac{1}{2}\lambda_1^2$

Since $\xi^{(2)}(1)$ is positive, we determine that

(3.17) $\quad\quad\quad \lambda_2 > \lambda_1(2-\lambda_1)$

and from (2.14) we see that $2 - \lambda_1 > 0$. We may therefore conclude that $\lambda_2$ is also positive.

Since $1 > \lambda_1$ we have $\lambda_1 > \lambda_1^2$ and from (3.17) we have

$$\lambda_2 > 2\lambda_1 - \lambda_1^2 > 2\lambda_1 - \lambda_1$$

and therefore we also have

(3.18) $\quad\quad\quad \lambda_2 > \lambda_1$

Continuing as above we have

$$\xi^{(3)}(s) = -\xi(s)\sum_{k=1}^{\infty}\frac{\sigma_k}{k}k(k-1)(k-2)s^{k-3} - 2\xi^{(1)}(s)\sum_{k=1}^{\infty}\frac{\sigma_k}{k}k(k-1)s^{k-2} - \xi^{(2)}(s)\sum_{k=1}^{\infty}\frac{\sigma_k}{k}ks^{k-1}$$

and using the identity $k(k-1)(k-2) = k(k+1)(k+2) - 6k(k+1) + 6k$ we obtain

$$\xi^{(3)}(s) = -\xi(s)\sum_{k=1}^{\infty}\frac{\sigma_k}{k}k(k+1)(k+2)s^{k-3} + 6\xi(s)\sum_{k=1}^{\infty}\frac{\sigma_k}{k}k(k+1)s^{k-3} - 6\xi(s)\sum_{k=1}^{\infty}\frac{\sigma_k}{k}ks^{k-3}$$

$$-2\xi^{(1)}(s)\sum_{k=1}^{\infty}\frac{\sigma_k}{k}k(k+1)s^{k-2} + 4\xi^{(1)}(s)\sum_{k=1}^{\infty}\frac{\sigma_k}{k}ks^{k-2} - \xi^{(2)}(s)\sum_{k=1}^{\infty}\frac{\sigma_k}{k}ks^{k-1}$$

In particular we have

$$\xi^{(3)}(1) = -\xi(1)\sum_{k=1}^{\infty}\frac{\sigma_k}{k}k(k+1)(k+2) + 6\xi(1)\sum_{k=1}^{\infty}\frac{\sigma_k}{k}k(k+1) - 6\xi(1)\sum_{k=1}^{\infty}\frac{\sigma_k}{k}k$$

$$-2\xi^{(1)}(1)\sum_{k=1}^{\infty}\frac{\sigma_k}{k}k(k+1) + 4\xi^{(1)}(1)\sum_{k=1}^{\infty}\frac{\sigma_k}{k}k - \xi^{(2)}(1)\sum_{k=1}^{\infty}\frac{\sigma_k}{k}k$$

Employing (3.16) this becomes

$$\xi^{(3)}(1) = -\frac{1}{2}\sum_{k=1}^{\infty}\frac{\sigma_k}{k}k(k+1)(k+2) + 3\sum_{k=1}^{\infty}\frac{\sigma_k}{k}k(k+1) - 3\sum_{k=1}^{\infty}\frac{\sigma_k}{k}k$$

$$-\lambda_1\sum_{k=1}^{\infty}\frac{\sigma_k}{k}k(k+1) + 2\lambda_1\sum_{k=1}^{\infty}\frac{\sigma_k}{k}k - \frac{1}{2}(\lambda_2 - 2\lambda_1 + \lambda_1^2)\sum_{k=1}^{\infty}\frac{\sigma_k}{k}k$$

and we obtain

$$\xi^{(3)}(1) = \lambda_3 - 3\lambda_2 + 3\lambda_1 + \lambda_1\lambda_2 - 2\lambda_1^2 + \frac{1}{2}\lambda_1(\lambda_2 - 2\lambda_1 + \lambda_1^2)$$

Since $\xi^{(3)}(1)$ is positive we have

$$\lambda_3 > 3(\lambda_2 - \lambda_1) - 3\lambda_1\left(\frac{1}{2}\lambda_2 - \lambda_1\right) - \frac{1}{2}\lambda_1^3$$

$$= 3(\lambda_2 - \lambda_1) - 3\lambda_1(\lambda_2 - \lambda_1) + \frac{3}{2}\lambda_1\lambda_2 - \frac{1}{2}\lambda_1^3$$

$$= 3(\lambda_2 - \lambda_1)(1-\lambda_1) + \frac{3}{2}\lambda_1\lambda_2 - \frac{1}{2}\lambda_1^3$$

$$= 3(\lambda_2 - \lambda_1)(1-\lambda_1) + \frac{1}{2}\lambda_1(3\lambda_2 - \lambda_1^2)$$

$$> 3(\lambda_2 - \lambda_1)(1-\lambda_1) + \frac{1}{2}\lambda_1(3\lambda_2 - \lambda_1) > 0$$

and we therefore deduce that $\lambda_3$ is positive. It remains to be determined whether or not the above inequality also implies that $\lambda_3 > \lambda_2$.

Using the other version of (1.7) gives us

(3.19) $$\xi^{(2)}(s) = -\xi(s)\sum_{k=1}^{\infty}(-1)^k(k-1)\sigma_k(s-1)^{k-2} - \xi^{(1)}(s)\sum_{k=1}^{\infty}(-1)^k\sigma_k(s-1)^{k-1}$$

and for $s=1$ we get

$$\xi^{(2)}(1) = -\xi(1)\sigma_2 + \xi^{(1)}(1)\sigma_1 = \frac{1}{2}[\sigma_1^2 - \sigma_2]$$

Hence we see that

(3.20) $\quad \sigma_1^2 > \sigma_2$

Letting $s=0$ in (3.19) results in

$$\xi^{(2)}(0) = -\frac{1}{2}\sum_{k=1}^{\infty}\frac{\sigma_k}{k}k(k-1) - \xi^{(1)}(s)\sum_{k=1}^{\infty}\frac{\sigma_k}{k}k = \frac{1}{2}\lambda_2 - \lambda_1 + \frac{1}{2}\lambda_1^2$$

which is the same as (3.16) since

$$\xi^{(n)}(s) = (-1)^n \xi^{(n)}(1-s)$$

and therefore

$$\xi^{(n)}(0) = (-1)^n \xi^{(n)}(1)$$

A further differentiation gives us

$$\xi^{(3)}(s) = -\xi(s)\sum_{k=1}^{\infty}(-1)^k \frac{\sigma_k}{k} k(k+1)(k+2)(s-1)^{k-3} + 6\xi(s)\sum_{k=1}^{\infty}(-1)^k \frac{\sigma_k}{k} k(k+1)(s-1)^{k-3}$$

$$-6\xi(s)\sum_{k=1}^{\infty}(-1)^k \frac{\sigma_k}{k} k (s-1)^{k-3} - 2\xi^{(1)}(s)\sum_{k=1}^{\infty}(-1)^k \frac{\sigma_k}{k} k(k+1)(s-1)^{k-2}$$

$$+4\xi^{(1)}(s)\sum_{k=1}^{\infty}(-1)^k \frac{\sigma_k}{k} k (s-1)^{k-2} - \xi^{(2)}(s)\sum_{k=1}^{\infty}(-1)^k \frac{\sigma_k}{k} k(s-1)^{k-1}$$

It will be seen in (6.2) below that

$$\xi^{(n+1)}(1) = \frac{1}{2}(-1)^n n! \sigma_{n+1} + \sum_{k=1}^{n}\binom{n}{k}(-1)^{n-k+1}\xi^{(k)}(1)(n-k+1)!\sigma_{n-k}$$

□

We now attempt to generalise the above specific results which we obtained for $\lambda_1, \lambda_2$ and $\lambda_3$.

It is shown in [22] that

(3.22) $$\frac{d^m}{dx^m} e^{f(x)} = e^{f(x)} Y_m\left(f^{(1)}(x), f^{(2)}(x), ..., f^{(m)}(x)\right)$$

where the (exponential) complete Bell polynomials $Y_n(x_1, ..., x_n)$ are defined by $Y_0 = 1$ and for $n \geq 1$

(3.23) $$Y_n(x_1, ..., x_n) = \sum_{\pi(n)} \frac{n!}{k_1! k_2! ... k_n!} \left(\frac{x_1}{1!}\right)^{k_1} \left(\frac{x_2}{2!}\right)^{k_2} ... \left(\frac{x_n}{n!}\right)^{k_n}$$

where the sum is taken over all partitions $\pi(n)$ of $n$, i.e. over all sets of integers $k_j$ such that

$$k_1 + 2k_2 + 3k_3 + ... + nk_n = n$$

The complete Bell polynomials have integer coefficients and the first five are set out below (Comtet [14, p.307])

(3.24) $$Y_1(x_1) = x_1$$

$$Y_2(x_1, x_2) = x_1^2 + x_2$$

$$Y_3(x_1, x_2, x_3) = x_1^3 + 3x_1 x_2 + x_3$$

$$Y_4(x_1, x_2, x_3, x_4) = x_1^4 + 6x_1^2 x_2 + 4x_1 x_3 + 3x_2^2 + x_4$$

$$Y_5(x_1, x_2, x_3, x_4, x_5) = x_1^5 + 10x_1^3 x_2 + 10x_1^2 x_3 + 15x_1 x_2^2 + 5x_1 x_4 + 10x_2 x_3 + x_5$$

Suppose that $h'(x) = h(x)g(x)$ and let $f(x) = \log h(x)$. We see that

$$f'(x) = \frac{h'(x)}{h(x)} = g(x)$$

and then using (3.22) above we have

$$\frac{d^m}{dx^m} h(x) = \frac{d^m}{dx^m} e^{\log h(x)} = h(x) Y_m\left(g(x), g^{(1)}(x), ..., g^{(m-1)}(x)\right)$$

We now write (1.7) as

$$\xi'(s) = \xi(s)g(s)$$

where $g(s) = -\sum_{k=1}^{\infty} \frac{\sigma_k}{k} k s^{k-1}$ and $g(1) = \lambda_1$. We then have

(3.25) $$\xi^{(m)}(s) = \xi(s) Y_m\left(g(s), g^{(1)}(s), ..., g^{(m-1)}(s)\right)$$

where

$$g^{(r)}(s) = -\sum_{k=1}^{\infty} \frac{\sigma_k}{k} k(k-1)\cdots(k-r) s^{k-1-r}$$

and

$$g^{(r)}(1) = -\sum_{k=1}^{\infty} \frac{\sigma_k}{k} k(k-1)\cdots(k-r)$$

Using (3.9)

$$k(k-1)...(k-r) = (-1)^{r+1} \sum_{j=1}^{r+1} \frac{(r+1)!}{j!} \binom{r}{j-1} (-1)^j k(k+1)...(k+j-1)$$

we have

$$g^{(r)}(1) = -(-1)^{r+1} \sum_{j=1}^{r+1} \frac{(r+1)!}{j!} \binom{r}{j-1} (-1)^j \sum_{k=1}^{\infty} \frac{\sigma_k}{k} k(k+1)\ldots(k+j-1)$$

and using (3.4) this becomes

$$= (-1)^{r+1} \sum_{j=1}^{r+1} \frac{(r+1)!}{j!} \binom{r}{j-1} (-1)^j (j-1)! \lambda_j$$

$$= (-1)^{r+1} (r+1)! \sum_{j=1}^{r+1} \binom{r}{j-1} (-1)^j \frac{\lambda_j}{j}$$

We then have

(3.26) $\qquad g^{(r)}(1) = (-1)^{r+1} r! \sum_{j=1}^{r+1} \binom{r+1}{j} (-1)^j \lambda_j$

and applying the binomial inversion formula

(3.27) $\qquad a_n = \sum_{k=0}^{n} \binom{n}{k} (-1)^k b_k \quad \Leftrightarrow \quad b_n = \sum_{k=0}^{n} \binom{n}{k} (-1)^k a_k$

to the following expression (where we define $\lambda_0 = g^{(-1)}(1) = 0$)

$$\frac{(-1)^r g^{(r-1)}(1)}{(r-1)!} = \sum_{j=0}^{r} \binom{r}{j} (-1)^j \lambda_j$$

we obtain

(3.28) $\qquad \lambda_j = \sum_{j=1}^{r} \binom{r}{j} \frac{g^{(j-1)}(1)}{(j-1)!}$

We have

$$g^{(j-1)}(s) = -\sum_{k=1}^{\infty} (-1)^k \frac{\sigma_k}{k} k(k-1) \cdots (k-j+1)(s-1)^{k-j}$$

which results in

$$g^{(j-1)}(1) = -(j-1)!(-1)^j \sigma_j$$

and therefore from (3.28) we obtain the well-known result [6]

(3.29) $$\lambda_j = -\sum_{j=1}^{r}(-1)^j \binom{r}{j}\sigma_j$$

Particular values of (3.26) are

$$g(1) = \lambda_1$$

$$g^{(1)}(1) = \lambda_2 - 2\lambda_1$$

$$g^{(2)}(1) = 6\lambda_1 - 6\lambda_2 + 2\lambda_3$$

Using (3.25) we see that

$$\xi^{(m)}(1) = \frac{1}{2}Y_m\left(g(1), g^{(1)}(1), \ldots, g^{(m-1)}(1)\right)$$

$$\xi^{(1)}(1) = \frac{1}{2}Y_1(g(1)) = \frac{1}{2}\lambda_1$$

$$\xi^{(2)}(1) = \frac{1}{2}Y_2\left(g(1), g^{(1)}(1)\right)$$

$$= \frac{1}{2}[g^2(1) + g^{(1)}(1)]$$

$$= \frac{1}{2}[\lambda_1^2 + \lambda_2 - 2\lambda_1]$$

$$\xi^{(3)}(1) = \frac{1}{2}Y_3\left(g(1), g^{(1)}(1), g^{(2)}(1)\right)$$

and, since $Y_3(x_1, x_2, x_3) = x_1^3 + 3x_1x_2 + x_3$, this becomes

$$= \frac{1}{2}[g^3(1) + 3g(1)g^{(1)}(1) + g^{(2)}(1)]$$

$$\xi^{(3)}(1) = \frac{1}{2}[\lambda_1^3 + 3\lambda_1(\lambda_2 - 2\lambda_1) + 6\lambda_1 - 6\lambda_2 + 2\lambda_3]$$

Riordan [28] reports that

$$(3.30) \quad Y_{m+1}(x_1,...,x_{m+1}) = \sum_{k=0}^{m}\binom{m}{k}Y_{m-k}(x_1,...,x_{m-k})x_{k+1} = \sum_{k=0}^{m}\binom{m}{k}Y_k(x_1,...,x_k)x_{m-k+1}$$

and for example we have

$$Y_2(x_1,x_2) = Y_1(x_1)x_1 + Y_0 x_2 = x_1^2 + x_2$$

From (3.25) we see that

$$\xi^{(m+1)}(s) = \xi(s)Y_{m+1}\big(g(s), g^{(1)}(s),..., g^{(m)}(s)\big)$$

$$= \sum_{k=0}^{m}\binom{m}{k}\xi(s)Y_k\big(g(s), g^{(1)}(s),..., g^{(k-1)}(s)\big)g^{(m-k)}(s)$$

We therefore have

$$(3.31) \quad \xi^{(m+1)}(s) = \sum_{k=0}^{m}\binom{m}{k}\xi^{(k)}(s)g^{(m-k)}(s)$$

which simply corresponds with the Leibniz differentiation rule using $\xi'(s) = \xi(s)g(s)$.

This gives us

$$\xi^{(2)}(s) = \xi(s)g^{(1)}(s) + \xi^{(1)}(s)g(s)$$

$$\xi^{(2)}(1) = \frac{1}{2}[\lambda_2 - 2\lambda_1] + \frac{1}{2}\lambda_1^2$$

We therefore have

$$\frac{1}{2}[\lambda_2 - 2\lambda_1] + \frac{1}{2}\lambda_1^2 > 0$$

which corresponds with (3.16).

□

We note from (1.6) that

$$\log \xi(s) = \log \Gamma\left(1+\frac{s}{2}\right) - \frac{s}{2}\log \pi + \log[(s-1)\varsigma(s)]$$

and therefore

$$g(s) = \frac{\xi'(s)}{\xi(s)} = \psi\left(1+\frac{s}{2}\right) - \frac{1}{2}\log \pi + \frac{1}{s-1} + \frac{\varsigma'(s)}{\varsigma(s)}$$

We will note from (4.1) that

$$\frac{\varsigma'(s)}{\varsigma(s)} + \frac{1}{s-1} = -\sum_{k=0}^{\infty} \eta_k (s-1)^k$$

and we therefore obtain

(3.32) $$g(s) = \frac{\xi'(s)}{\xi(s)} = \frac{1}{2}\psi\left(1+\frac{s}{2}\right) - \frac{1}{2}\log \pi - \sum_{k=0}^{\infty} \eta_k (s-1)^k$$

and we see that

$$g(1) = g^{(0)}(1) = \frac{1}{2}\psi\left(\frac{3}{2}\right) - \frac{1}{2}\log \pi - \eta_0 = \lambda_1$$

For $r \geq 1$

$$g^{(r)}(s) = \frac{1}{2^{r+1}}\psi^{(r)}\left(1+\frac{s}{2}\right) - \sum_{k=0}^{\infty} \eta_k k(k-1)\cdots(k-r+1)(s-1)^{k-r}$$

$$g^{(r)}(1) = \frac{1}{2^{r+1}}\psi^{(r)}\left(\frac{3}{2}\right) - r!\eta_r$$

$$= \frac{1}{2^{r+1}}\psi^{(r)}\left(\frac{3}{2}\right) - r!\eta_r - \delta_{r,0}\frac{1}{2}\log \pi$$

where we have introduced the Kronecker delta $\delta_{r,0}$ to ensure that this equation is also valid for $r = 0$

$$g^{(r)}(1) = (-1)^{r+1} r! \sum_{j=1}^{r+1} \binom{r+1}{j} (-1)^j \lambda_j$$

This may be written as

$$(-1)^{r+1} r! \sum_{j=0}^{r+1} \binom{r+1}{j}(-1)^j \lambda_j = \frac{1}{2^{r+1}} \psi^{(r)}\left(\frac{3}{2}\right) - r!\eta_r - \delta_{r,0} \frac{1}{2}\log \pi$$

where we start the summation at $j = 0$ by defining $\lambda_0 = 0$.

We now let $r \to r-1$ and define $\eta_{-1} = 0$ and $\psi^{(-1)}\left(\frac{3}{2}\right) = 0$

$$(-1)^r (r-1)! \sum_{j=0}^{r} \binom{r}{j}(-1)^j \lambda_j = \frac{1}{2^r} \psi^{(r-1)}\left(\frac{3}{2}\right) - (r-1)!\eta_{r-1} - \delta_{r-1,0} \frac{1}{2}\log \pi$$

$$\sum_{j=0}^{r} \binom{r}{j}(-1)^j \lambda_j = (-1)^r \left[\frac{r}{2^r r!} \psi^{(r-1)}\left(\frac{3}{2}\right) - \eta_{r-1} - \delta_{r-1,0} \frac{r}{2r!}\log \pi \right]$$

and applying the binomial inversion formula (3.27) we obtain

$$\lambda_r = \sum_{j=0}^{r} \binom{r}{j}\left[\frac{j}{2^j j!} \psi^{(j-1)}\left(\frac{3}{2}\right) - \eta_{j-1} - \delta_{j-1,0} \frac{j}{2j!}\log \pi \right]$$

$$\lambda_r = \sum_{j=1}^{r} \binom{r}{j}\left[\frac{j}{2^j j!} \psi^{(j-1)}\left(\frac{3}{2}\right) - \eta_{j-1}\right] - \frac{1}{2}\binom{r}{1} \log \pi$$

(3.33) $$\lambda_r = \sum_{j=2}^{r} \binom{r}{j}\frac{j}{2^j j!} \psi^{(j-1)}\left(\frac{3}{2}\right) + \frac{1}{2}\binom{r}{j}\psi\left(\frac{3}{2}\right) - \sum_{j=1}^{r}\binom{r}{j}\eta_{j-1} - \frac{1}{2}\binom{r}{1}\log \pi$$

From (3.5.1) we have

$$\psi^{(j-1)}(3/2) = (-1)^j (j-1)!\left(\left[2^j - 1\right]\varsigma(j) - 2^j\right)$$

and we obtain

$$\lambda_r = \sum_{j=2}^{r}\binom{r}{j}(-1)^j \left[1 - \frac{1}{2^j}\right]\varsigma(j) - \sum_{j=2}^{r}\binom{r}{j}(-1)^j + \frac{1}{2}\binom{r}{1}\psi\left(\frac{3}{2}\right) - \sum_{j=1}^{r}\binom{r}{j}\eta_{j-1} - \frac{1}{2}\binom{r}{1}\log \pi$$

Using (2.11)

$$\psi\left(\frac{3}{2}\right) = -\gamma - 2\log 2 + 2$$

and noting that

$$\sum_{j=2}^{r}\binom{r}{j}(-1)^j = \sum_{j=0}^{r}\binom{r}{j}(-1)^j + r - 1 = r - 1$$

we simply obtain another derivation of Coffey's result [C1] for $r \geq 2$

(3.34) $$\lambda_r = \sum_{j=2}^{r}\binom{r}{j}(-1)^j\left[1-\frac{1}{2^j}\right]\varsigma(j) + 1 - \sum_{j=1}^{r}\binom{r}{j}\eta_{j-1} + 1 - \frac{1}{2}r(\gamma + 2\log 2 + \log \pi)$$

(and we note that the term +1 has been inadvertently omitted in equation (3) in [7] and in equation (12) in [13]).

A much simpler derivation of (3.34) is shown below.

We see that (3.33) may be written as

$$\lambda_r = \sum_{j=2}^{r}\binom{r}{j}\frac{1}{(j-1)!}\left[\frac{1}{2^j}\psi^{(j-1)}\left(\frac{3}{2}\right) - (j-1)!\eta_{j-1}\right] + r\gamma + \frac{1}{2}r\psi\left(\frac{3}{2}\right) - \frac{1}{2}r\log\pi$$

$$= \sum_{j=2}^{r}\binom{r}{j}\frac{1}{(j-1)!}\left[\frac{1}{2^j}\psi^{(j-1)}\left(\frac{3}{2}\right) - (j-1)!\eta_{j-1}\right] + \frac{1}{2}r(\gamma - 2\log 2 + 2 - \log\pi)$$

and we note from (3.32) that

$$\frac{d}{ds}\log \xi(s) + \frac{1}{2}\log \pi = \frac{1}{2}\psi\left(1+\frac{s}{2}\right) - \sum_{k=0}^{\infty}\eta_k(s-1)^k$$

Differentiating this $j-1$ times gives us

$$\frac{d^j}{ds^j}\log \xi(s) + \delta_{j,1}\frac{1}{2}\log \pi = \frac{1}{2^j}\psi^{(j-1)}\left(1+\frac{s}{2}\right) - \sum_{k=0}^{\infty}\eta_k k(k-1)\cdots(k-j)(s-1)^{k-j+1}$$

and evaluating this at $s=1$ results in

$$\frac{d^j}{ds^j}\log \xi(s)\bigg|_{s=1} + \delta_{j,1}\frac{1}{2}\log \pi = \frac{1}{2^j}\psi^{(j-1)}\left(\frac{3}{2}\right) - (j-1)!\eta_{j-1}$$

We then have

$$\sum_{j=2}^{r}\binom{r}{j}\frac{1}{(j-1)!}\left[\frac{1}{2^j}\psi^{(j-1)}\left(\frac{3}{2}\right) - (j-1)!\eta_{j-1}\right] + r\gamma + \frac{1}{2}r\psi\left(\frac{3}{2}\right) - \frac{1}{2}r\log\pi$$

$$= \sum_{j=2}^{r}\binom{r}{j}\frac{1}{(j-1)!}\frac{d^j}{ds^j}\log \xi(s)\bigg|_{s=1} -\frac{1}{2}r(\gamma +2\log 2+\log \pi)$$

Applying the Leibniz differentiation rule to (1.4) we see that

$$\lambda_r = \frac{1}{(r-1)!}\frac{d^r}{ds^r}[s^{r-1}\log \xi(s)]\bigg|_{s=1}$$

$$= \frac{1}{(r-1)!}\sum_{j=0}^{r}\binom{r}{j}\frac{d^j}{ds^j}\log \xi(s)\frac{d^{r-j}}{ds^{r-j}}s^{r-1}\bigg|_{s=1}$$

$$= \sum_{j=0}^{r}\binom{r}{j}\frac{1}{(j-1)!}\frac{d^j}{ds^j}\log \xi(s)\bigg|_{s=1}$$

$$= \sum_{j=2}^{r}\binom{r}{j}\frac{1}{(j-1)!}\frac{d^j}{ds^j}\log \xi(s)\bigg|_{s=1} +r\frac{d}{ds}\log \xi(s)\bigg|_{s=1}$$

$$= \sum_{j=2}^{r}\binom{r}{j}\frac{1}{(j-1)!}\left[\frac{1}{2^j}\psi^{(j-1)}\left(\frac{3}{2}\right)-(j-1)!\eta_{j-1}\right]+\frac{1}{2}r(\gamma -2\log 2+2-\log \pi)$$

which corresponds with (3.33) above. $\square$

We have

$$\lambda_r = \sum_{j=2}^{r}\binom{r}{j}\frac{1}{(j-1)!}\left[\frac{1}{2^j}\psi^{(j-1)}\left(\frac{3}{2}\right)-(j-1)!\eta_{j-1}\right]+\frac{1}{2}r(\gamma -2\log 2+2-\log \pi)$$

Using (3.5.2)

$$\psi^{(j-1)}(3/2)=(-1)^j(j-1)!\varsigma(j,3/2)$$

we may write this as

(3.35) $$\lambda_r = \sum_{j=2}^{r}\binom{r}{j}(-1)^j\left[\frac{\varsigma(j,3/2)}{2^j}-(-1)^j\eta_{j-1}\right]+\frac{1}{2}r(\gamma -2\log 2+2-\log \pi)$$

and it may be possible to simplify the first summation by using the following Hermite integral representation of the Hurwitz zeta function (see for example [31, p.120])

$$\varsigma(s,u)=\frac{u^{-s}}{2}+\frac{u^{1-s}}{s-1}+i\int_0^\infty \frac{(u+ix)^{-s}-(u-ix)^{-s}}{e^{2\pi x}-1}dx$$

It may be noted that

$$\gamma - 2\log 2 + 2 - \log \pi > 0$$

## 4. The eta constants

The eta constants $\eta_n$ are defined by reference to the logarithmic derivative of the Riemann zeta function

(4.1) $$\frac{d}{ds}[\log \varsigma(s)] = \frac{\varsigma'(s)}{\varsigma(s)} = -\frac{1}{s-1} - \sum_{k=0}^{\infty} \eta_k (s-1)^k \qquad |s-1|<3$$

and we may also note that this is equivalent to

(4.2) $$\frac{d}{ds}\log[(s-1)\varsigma(s)] = \frac{\varsigma'(s)}{\varsigma(s)} + \frac{1}{s-1} = -\sum_{k=0}^{\infty} \eta_k (s-1)^k$$

We then see that

$$\frac{d^{n+1}}{ds^{n+1}}\log[(s-1)\varsigma(s)] = -\sum_{k=0}^{\infty} \eta_k k(k-1)\ldots(k-n+1)(s-1)^{k-n}$$

and hence we get

(4.3) $$\lim_{s \to 1} \frac{d^{n+1}}{ds^{n+1}}\log[(s-1)\varsigma(s)] = -n!\eta_n$$

We see from (4.2) that

$$\frac{d}{ds}[(s-1)\varsigma(s)] = -\sum_{k=0}^{\infty} \eta_k (s-1)^k [(s-1)\varsigma(s)]$$

and thus

$$\frac{d^{n+1}}{ds^{n+1}}[(s-1)\varsigma(s)] = -\frac{d^n}{ds^n}\sum_{k=0}^{\infty} \eta_k (s-1)^k [(s-1)\varsigma(s)]$$

We consider the value at $s=1$

$$\left.\frac{d^{n+1}}{ds^{n+1}}[(s-1)\varsigma(s)]\right|_{s=1} = -\lim_{s \to 1} \frac{d^n}{ds^n}\sum_{k=0}^{\infty} \eta_k (s-1)^k [(s-1)\varsigma(s)]$$

and using (2.10)

$$\left. \frac{d^{n+1}}{ds^{n+1}}[(s-1)\varsigma(s)] \right|_{s=1} = (-1)^n (n+1)\gamma_n$$

we obtain via the Leibniz differentiation formula

$$= -n!\eta_n + \sum_{j=1}^{n} \binom{n}{j}(n-j)!\eta_{n-j}(-1)^j j\gamma_{j-1}$$

where we have isolated the first term.

$$= -n!\eta_n + n!\sum_{j=1}^{n} \frac{(-1)^j j\eta_{n-j}\gamma_{j-1}}{j!}$$

Hence we obtain the recurrence relation

(4.4) $$(-1)^n (n+1)\gamma_n = -n!\eta_n + n!\sum_{j=1}^{n} \frac{(-1)^j j\eta_{n-j}\gamma_{j-1}}{j!}$$

where for $n = 0, 1, 2$ we have respectively:

$$\gamma = -\eta_0$$

$$2\gamma_1 = \eta_1 + \eta_0 \gamma = \eta_1 - \gamma^2$$

$$\frac{3}{2}\gamma_2 = -\eta_2 - 3\gamma\gamma_1 - \gamma^3$$

Equation (4.4) is equivalent to Coffey's recurrence relation (4.5) below.

Letting $u = 1$ in (2.1) gives us

$$(s-1)\varsigma(s) = 1 + \sum_{n=0}^{\infty} \frac{(-1)^n}{n!}\gamma_n(s-1)^{n+1}$$

and defining

$$g(s) = (s-1)\varsigma(s)$$

$$L(s) = \log[(s-1)\varsigma(s)]$$

we have

$$\frac{d}{ds}(s-1)\varsigma(s) = \sum_{n=0}^{\infty} \frac{(-1)^n}{n!}\gamma_n(n+1)(s-1)^n$$

and thus

$$g^{(1)}(1) = \gamma_0 = \gamma$$

$$g^{(k+1)}(1) = (-1)^k (k+1)\gamma_k$$

We see that

$$\frac{d}{ds}\log[(s-1)\varsigma(s)] = \frac{(s-1)\varsigma'(s) + \varsigma(s)}{(s-1)\varsigma(s)}$$

which may be written as

$$(s-1)\varsigma(s)L^{(1)}(s) = (s-1)\varsigma'(s) + \varsigma(s) = \frac{d}{ds}[(s-1)\varsigma(s)]$$

or

$$g(s)L^{(1)}(s) = g^{(1)}(s)$$

We then have

$$\frac{d^{(n)}}{ds^{(n)}}\left[g(s)L^{(1)}(s)\right] = g^{(n+1)}(s)$$

and applying the Leibniz rule gives us

$$g^{(n+1)}(s) = \sum_{k=0}^{n}\binom{n}{k}g^{(k)}(s)L^{(n+1-k)}(s)$$

Letting $s = 1$ gives us

$$(-1)^n(n+1)\gamma_n = g^{(0)}(1)L^{(n+1)}(1) - \sum_{k=1}^{n}\binom{n}{k}(-1)^{k-1}k\gamma_{k-1}(n-k)!\eta_{n-k}$$

$$= -n!\eta_n - \sum_{k=1}^{n}\binom{n}{k}(-1)^k k\gamma_{k-1}(n-k)!\eta_{n-k}$$

$$= -n!\eta_n - n!\sum_{k=1}^{n}\frac{(-1)^k}{(k-1)!}\gamma_{k-1}\eta_{n-k}$$

Hence we have

$$n!\eta_n = (-1)^{n+1}(n+1)\gamma_n - n!\sum_{k=1}^{n}\frac{(-1)^k}{(k-1)!}\gamma_{k-1}\eta_{n-k}$$

and reindexing the sum gives us

(4.5) $$n!\eta_n = (-1)^{n+1}(n+1)\gamma_n + (-1)^{n+1}n!\sum_{k=0}^{n-1}\frac{(-1)^{k+1}}{(n-k-1)!}\gamma_{n-k-1}\eta_k$$

as originally reported by Coffey [5].

A formula for the Stieltjes constants in terms of the (exponential) complete Bell polynomials containing the eta constants $\eta_n$ as the arguments is shown below in equation (6.1).

$\square$

Eliminating $\log\xi(s)$ from (1.6) and (1.7) results in

$$\log\Gamma\left(1+\frac{s}{2}\right) - \frac{s}{2}\log\pi + \log[(s-1)\varsigma(s)] = -\log 2 - \sum_{k=1}^{\infty}(-1)^k\frac{\sigma_k}{k}(s-1)^k$$

and, since $\xi(s) = \xi(1-s)$, we have the equivalent representation

$$\log\Gamma\left(1+\frac{s}{2}\right) - \frac{s}{2}\log\pi + \log[(s-1)\varsigma(s)] = -\log 2 - \sum_{k=1}^{\infty}\frac{\sigma_k}{k}s^k$$

For convenience, we define $L(s)$ by $L(s) = \log[(s-1)\varsigma(s)]$

Upon differentiating the above two equations we obtain

$$\frac{1}{2}\psi\left(1+\frac{s}{2}\right) - \frac{1}{2}\log\pi + L^{(1)}(s) = -\sum_{k=1}^{\infty}(-1)^k\sigma_k(s-1)^{k-1}$$

$$\frac{1}{2}\psi\left(1+\frac{s}{2}\right) - \frac{1}{2}\log\pi + L^{(1)}(s) = -\sum_{k=1}^{\infty}\sigma_k s^{k-1}$$

For $n \geq 1$ we have the higher derivatives

$$\frac{1}{2^{n+1}}\psi^{(n)}\left(1+\frac{s}{2}\right) + L^{(n+1)}(s) = -\sum_{k=1}^{\infty}(-1)^k\sigma_k(k-1)(k-2)\ldots(k-n)(s-1)^{k-1-n}$$

$$\frac{1}{2^{n+1}}\psi^{(n)}\left(1+\frac{s}{2}\right) + L^{(n+1)}(s) = -\sum_{k=1}^{\infty}\sigma_k(k-1)(k-2)\ldots(k-n)s^{k-1-n}$$

and therefore with $s = 1$ we obtain

$$\frac{1}{2^{n+1}}\psi^{(n)}(3/2)+L^{(n+1)}(1)=-(-1)^{n+1}n!\sigma_{n+1}$$

We have

$$\psi^{(n)}(3/2)=(-1)^{n+1}n![2^{n+1}-1]\varsigma(n+1)+(-1)^n n!2^{n+1}$$

and hence we get

$$(-1)^{n+1}n!\left[1-\frac{1}{2^{n+1}}\right]\varsigma(n+1)+(-1)^n n!+L^{(n+1)}(1)=-(-1)^{n+1}n!\sigma_{n+1}$$

Using $L^{(n+1)}(1)=-n!\eta_n$ results in

$$\sigma_{n+1}=(-1)^{n+1}\eta_n-\left[1-\frac{1}{2^{n+1}}\right]\varsigma(n+1)+1$$

Letting $s=0$ gives us

$$\frac{1}{2^{n+1}}\psi^{(n)}(1)+L^{(n+1)}(0)=-n!\sigma_{n+1}$$

and therefore we get

$$L^{(n+1)}(0)=-n!\sigma_{n+1}-\frac{1}{2^{n+1}}(-1)^{n+1}n!\varsigma(n+1)$$

and substituting for $\sigma_{n+1}$ gives us

$$=-n!\left((-1)^{n+1}\eta_n-\left[1-\frac{1}{2^{n+1}}\right]\varsigma(n+1)+1\right)-\frac{1}{2^{n+1}}(-1)^{n+1}n!\varsigma(n+1)$$

Hence we get for $n\geq 1$

$$L^{(n+1)}(0)=(-1)^n n!\eta_n+\left(1-\frac{1}{2^{n+1}}\left[1-(-1)^n\right]\right)n!\varsigma(n+1)-n!$$

and with $n=1$ we have

$$L^{(2)}(0)=\frac{1}{2}\varsigma(2)-\eta_1-1$$

We will also see from (4.6) below that

$$L^{(2)}(0) = -2\varsigma''(0) - \log^2(2\pi) - 1$$

and we therefore obtain

$$\varsigma''(0) = \frac{1}{2}\eta_1 - \frac{1}{4}\varsigma(2) - \frac{1}{2}\log^2(2\pi)$$

Substituting (4.4) $\eta_1 = 2\gamma_1 + \gamma^2$ we then have an expression for $\varsigma''(0)$

$$\varsigma''(0) = \gamma_1 + \frac{1}{2}\gamma^2 - \frac{1}{4}\varsigma(2) - \frac{1}{2}\log^2(2\pi)$$

as previously derived by Ramanujan [3] and Apostol [1].

We have

$$L^{(1)}(s) = -\sum_{k=0}^{\infty} \eta_k (s-1)^k$$

and hence

$$L^{(1)}(1) = -\eta_0 = \gamma$$

$$L^{(n+1)}(s) = -\sum_{k=0}^{\infty} \eta_k k(k-1)...(k-n+1)(s-1)^{k-n}$$

Coffey [13] has shown that the sequence $(\eta_n)$ has strict sign alteration

$$\eta_n = (-1)^{n+1} \varepsilon_n$$

where $\varepsilon_n$ are positive constants and therefore we have

$$L^{(n+1)}(1) = -n!\eta_n = (-1)^n n!\varepsilon_n$$

We therefore note that $L^{(1)}(1)$ is positive and $L^{(2)}(1)$ is negative and the signs strictly alternate thereafter.

We see that

$$L^{(1)}(0) = -\sum_{k=0}^{\infty} (-1)^k \eta_k = \sum_{k=0}^{\infty} \varepsilon_k \text{ is positive}$$

In fact we have

$$L^{(1)}(0) = -\sum_{k=0}^{\infty}(-1)^k \eta_k = \log(2\pi) - 1$$

$$L^{(n+1)}(0) = -\sum_{k=0}^{\infty}\eta_k k(k-1)\ldots(k-n+1)(-1)^{k-n} = (-1)^n \sum_{k=0}^{\infty}\varepsilon_k k(k-1)\ldots(k-n+1)$$

and therefore we note that the signs of $L^{(n+1)}(0)$ also strictly alternate.

We see that

$$L^{(2)}(s) = \frac{(s-1)\varsigma''(s) + 2\varsigma'(s)}{(s-1)\varsigma(s)} - \left[L^{(1)}(s)\right]^2$$

Referring to (2.9) we then obtain

$$L^{(2)}(1) = -2\gamma_1 - \left[L^{(1)}(1)\right]^2 = -2\gamma_1 - \gamma^2$$

and we therefore obtain (as in 4.4))

$$\eta_1 = 2\gamma_1 + \gamma^2$$

Since $\eta_1 \geq 0$ we deduce that

$$2\gamma_1 + \gamma^2 \geq 0$$

Since $\eta_2 \leq 0$ we see that

$$-\frac{3}{2}\gamma_2 - \gamma^3 - 3\gamma\gamma_1 \leq 0$$

We have

$$L^{(2)}(0) = -2\varsigma^{(2)}(0) + 4\varsigma^{(1)}(0) - \left[L^{(1)}(0)\right]^2$$

$$= -2\varsigma^{(2)}(0) + 4\varsigma^{(1)}(0) - \left[\log(2\pi) - 1\right]^2$$

$$= -2\varsigma^{(2)}(0) - 2\log(2\pi) - \left[\log(2\pi) - 1\right]^2$$

and therefore

(4.6) $\quad L^{(2)}(0) = -2\varsigma^{(2)}(0) - \log^2(2\pi) - 1$

## 5. A formula for the Stieltjes constants in terms of the higher derivatives of the Riemann zeta function $\varsigma^{(n)}(0)$.

We recall Hasse's formula (2.4) with $u = 1$

$$(s-1)\varsigma(s) = (s-1)\varsigma(s,1) = \sum_{n=0}^{\infty} \frac{1}{n+1} \sum_{k=0}^{n} \binom{n}{k} \frac{(-1)^k}{(1+k)^{s-1}}$$

and with $s \to s - 1$ we have

$$s\varsigma(1-s) = -\sum_{n=0}^{\infty} \frac{1}{n+1} \sum_{k=0}^{n} \binom{n}{k} (-1)^k (1+k)^s$$

We refer to the functional equation for the Riemann zeta function (1.2)

$$2(2\pi)^{-s} \Gamma(s) \cos(\pi s / 2) \varsigma(s) = \varsigma(1-s)$$

and, multiplying by $s$, we see that

$$f(s) = s\varsigma(1-s) = 2(2\pi)^{-s} \Gamma(s+1) \cos(\pi s / 2) \varsigma(s) = -\sum_{n=0}^{\infty} \frac{1}{n+1} \sum_{k=0}^{n} \binom{n}{k} (-1)^k (1+k)^s$$

Differentiation results in

$$f^{(p)}(s) = -\sum_{n=0}^{\infty} \frac{1}{n+1} \sum_{k=0}^{n} \binom{n}{k} (-1)^k (1+k)^s \log^p(1+k)$$

and we have the particular value at $s = 0$

$$f^{(p)}(0) = -\sum_{n=0}^{\infty} \frac{1}{n+1} \sum_{k=0}^{n} \binom{n}{k} (-1)^k \log^p(1+k)$$

We recall the following expression for the Stieltjes constants (2.8)

$$\gamma_p(u) = -\frac{1}{p+1} \sum_{n=0}^{\infty} \frac{1}{n+1} \sum_{k=0}^{n} \binom{n}{k} (-1)^k \log^{p+1}(u+k)$$

which gives us

$$\gamma_p = \gamma_p(1) = -\frac{1}{p+1}\sum_{n=0}^{\infty}\frac{1}{n+1}\sum_{k=0}^{n}\binom{n}{k}(-1)^k \log^{p+1}(1+k)$$

and we therefore obtain

(5.1) $\qquad f^{(p+1)}(0) = (p+1)\gamma_p$

We have

$$\log f(s) = \log 2 - s\log(2\pi) + \log\Gamma(s+1) + \log\cos(\pi s/2) + \log\varsigma(s)$$

and differentiation results in

$$\frac{f'(s)}{f(s)} = -\log(2\pi) + \psi(s+1) - \frac{\pi}{2}\tan(\pi s/2) + \frac{\varsigma'(s)}{\varsigma(s)}$$

We note that

$$f(0) = 2\varsigma(0) = -1$$

$$f'(0) = f(0)\left[-\log(2\pi) + \psi(1) + \frac{\varsigma'(0)}{\varsigma(0)}\right] = \log(2\pi) + \gamma + 2\varsigma'(0)$$

and, referring to (5.1), we see that

$$\log(2\pi) + \gamma + 2\varsigma'(0) = \gamma_0 = \gamma$$

We then deduce the well-known result

(5.2) $\qquad \varsigma'(0) = -\frac{1}{2}\log(2\pi)$

We may write

$$f'(s) = g(s)f(s)$$

where

$$g(s) = \frac{d}{ds}\left[\log(2\pi)^{-s} + \log\Gamma(s+1) + \log\cos(\pi s/2) + \log\varsigma(s)\right]$$

$$= -\log(2\pi) + \psi(s+1) - \frac{\pi}{2}\tan(\pi s/2) + \frac{\varsigma'(s)}{\varsigma(s)}$$

We note that

$$g(0) = -\gamma \text{ and } f'(0) = \gamma$$

$$f'(s) = \left[-\log(2\pi) + \psi(s+1) - \frac{\pi}{2}\tan(\pi s/2) + \frac{\varsigma'(s)}{\varsigma(s)}\right]f(s)$$

$$f''(s) = \left[-\log(2\pi) + \psi(s+1) - \frac{\pi}{2}\tan(\pi s/2) + \frac{\varsigma'(s)}{\varsigma(s)}\right]f'(s)$$

$$+ \left[\psi'(s+1) - \left(\frac{\pi}{2}\right)^2 \sec^2(\pi s/2) + \frac{\varsigma(s)\varsigma''(s) - [\varsigma'(s)]^2}{\varsigma^2(s)}\right]f(s)$$

We then see that

$$f''(0) = -\gamma^2 - \left[\psi'(1) - \left(\frac{\pi}{2}\right)^2 + \frac{\varsigma(0)\varsigma''(0) - [\varsigma'(0)]^2}{\varsigma^2(0)}\right]$$

We know from (3.5.2) that

$$\psi'(1) = \varsigma(2) = \frac{\pi^2}{6}$$

and we then have

$$f''(0) = -\gamma^2 + \frac{\pi^2}{12} + 2\varsigma''(0) + \log^2(2\pi)$$

We also have from (5.1)

$$f''(0) = 2\gamma_1$$

and we thus obtain

(5.3) $$\varsigma''(0) = \gamma_1 + \frac{1}{2}\gamma^2 - \frac{1}{24}\pi^2 - \frac{1}{2}\log^2(2\pi)$$

This concurs with the result previously obtained by Ramanujan [3] and Apostol [1].

As noted by Coffey [9] for differentiable functions $f(s)$ and $g(s)$ such that

$$f'(s) = f(s)g(s)$$

then we have in terms of the (exponential) complete Bell polynomials

$$f^{(n)}(s) = f(s)Y_n\left(g(s), g'(s), \ldots, g^{(n-1)}(s)\right)$$

where, in the case under review, we have

$$g(s) = \frac{d}{ds}\left[\log(2\pi)^{-s} + \log\Gamma(s+1) + \log\cos(\pi s/2) + \log\varsigma(s)\right] = \frac{d}{ds}\sum_{p=1}^{4}\log g_p(s)$$

$$f^{(n)}(s) = f(s)Y_n\left(\frac{d}{ds}\sum_{p=1}^{4}\log g_p(s), \frac{d^2}{ds^2}\sum_{p=1}^{4}\log g_p(s), \ldots, \frac{d^{n-1}}{ds^{n-1}}\sum_{p=1}^{4}\log g_p(s)\right)$$

We have [28]

(5.4) $$Y_n(x_1 + y_1, \ldots, x_n + y_n) = \sum_{k=0}^{n}\binom{n}{k}Y_{n-k}(x_1, \ldots, x_{n-k})Y_k(y_1, \ldots, y_k)$$

and this may be generalised to Bell polynomials with additional arguments as follows

$$Y_n(a_1 + b_1 + c_1 + d_1, \ldots, a_n + b_n + c_n + d_n) = \sum_{k=0}^{n}\binom{n}{k}Y_{n-k}(a_1, \ldots, a_{n-k})Y_k(b_1 + c_1 + d_1, \ldots, b_k + c_k + d_k)$$

$$Y_k(b_1 + c_1 + d_1, \ldots, b_k + c_k + d_k) = \sum_{j=0}^{k}\binom{k}{j}Y_{k-j}(b_1, \ldots, b_{k-j})Y_j(c_1 + d_1, \ldots, c_j + d_j)$$

$$Y_j(c_1 + d_1, \ldots, c_j + d_j) = \sum_{l=0}^{j}\binom{j}{l}Y_{j-l}(c_1, \ldots, c_{j-l})Y_l(d_1, \ldots, d_l)$$

and we end up with

$$Y_n(a_1 + b_1 + c_1 + d_1, \ldots, a_n + b_n + c_n + d_n)$$

$$= \sum_{k=0}^{n}\binom{n}{k}Y_{n-k}(a_1, \ldots, a_{n-k})\sum_{j=0}^{k}\binom{k}{j}Y_{k-j}(b_1, \ldots, b_{k-j})\sum_{l=0}^{j}\binom{j}{l}Y_{j-l}(c_1, \ldots, c_{j-l})Y_l(d_1, \ldots, d_l)$$

We note from (A.5) in Appendix A that

$$\frac{d^m}{dx^m}e^{h(x)} = e^{h(x)}Y_m\left(h^{(1)}(x), h^{(2)}(x), \ldots, h^{(m)}(x)\right)$$

and letting $h(x) \to \log h(x)$ we have

$$\frac{d^m}{dx^m}h(x) = \frac{d^m}{dx^m}e^{\log h(x)} = e^{\log h(x)}Y_m\left(\frac{d}{dx}\log h(x), \frac{d^2}{dx^2}\log h(x), \ldots, \frac{d^m}{dx^m}\log h(x)\right)$$

and we see that

$$\frac{1}{h(x)}\frac{d^m}{dx^m}h(x) = Y_m\left(\frac{d}{dx}\log h(x), \frac{d^2}{dx^2}\log h(x), \ldots, \frac{d^m}{dx^m}\log h(x)\right)$$

We therefore have

$$f^{(n)}(s) = f(s)Y_n\left(\frac{d}{ds}\sum_{p=1}^{4}\log g_p(s), \frac{d^2}{ds^2}\sum_{p=1}^{4}\log g_p(s), \ldots, \frac{d^{n-1}}{ds^{n-1}}\sum_{p=1}^{4}\log g_p(s)\right)$$

$$= \frac{1}{g_1(s)g_2(s)g_3(s)g_4(s)}\sum_{k=0}^{n}\binom{n}{k}\frac{d^{n-k}}{ds^{n-k}}g_1(s)\sum_{j=0}^{k}\binom{k}{j}\frac{d^{k-j}}{ds^{k-j}}g_2(s)\sum_{l=0}^{j}\binom{j}{l}\frac{d^{j-l}}{ds^{j-l}}g_3(s)\frac{d^l}{ds^l}g_4(s)$$

and with $s = 0$ this becomes using (5.1)

$$\gamma_{n-1} = 2\sum_{k=0}^{n}\binom{n}{k}\frac{d^{n-k}}{ds^{n-k}}g_1(0)\sum_{j=0}^{k}\binom{k}{j}\frac{d^{k-j}}{ds^{k-j}}g_2(0)\sum_{l=0}^{j}\binom{j}{l}\frac{d^{j-l}}{ds^{j-l}}g_3(0)\frac{d^l}{ds^l}g_4(0)$$

where

$$g(s) = \frac{d}{ds}\left[\log(2\pi)^{-s} + \log\Gamma(s+1) + \log\cos(\pi s/2) + \log\varsigma(s)\right] = \frac{d}{ds}\sum_{p=1}^{4}\log g_p(s)$$

Having regard to the four components $g_p(s)$, the following derivatives are easily computed

$$\frac{d^m}{ds^m}(2\pi)^{-s} = (-1)^m(2\pi)^{-s}\log^m(2\pi)$$

$$\frac{d^m}{ds^m}(2\pi)^{-s}\bigg|_{s=0} = (-1)^m\log^m(2\pi)$$

$$\frac{d^m}{ds^m}\cos\left(\frac{\pi s}{2}\right) = \left(\frac{\pi}{2}\right)^m \cos\left(\frac{\pi s}{2} + \frac{m\pi}{2}\right)$$

$$\left.\frac{d^m}{ds^m}\cos\left(\frac{\pi s}{2}\right)\right|_{s=0} = \left(\frac{\pi}{2}\right)^m \cos\left(\frac{m\pi}{2}\right)$$

and we obtain

(5.5)
$$\gamma_{n-1} = 2\sum_{k=0}^{n}\binom{n}{k}\Gamma^{(n-k)}(1)\sum_{j=0}^{k}\binom{k}{j}(-1)^{k-j}\log^{k-j}(2\pi)\sum_{l=0}^{j}\binom{j}{l}\left(\frac{\pi}{2}\right)^{j-l}\cos\left(\frac{(j-l)\pi}{2}\right)\varsigma^{(l)}(0)$$

This extends the formula previously obtained by Apostol [1]. We note from (A.7) in Appendix A that the derivatives of the gamma function may be expressed in terms of the (exponential) complete Bell polynomials

$$\Gamma^{(m)}(1) = Y_m(-\gamma, x_1, \ldots, x_{m-1})$$

where $x_p = (-1)^{p+1} p!\varsigma(p+1)$.

$\square$

Using (3.3) we see that

$$\frac{d}{ds}[(s-1)\varsigma(s)] = -(s-1)\varsigma(s)\frac{d}{ds}\left[\sum_{k=1}^{\infty}(-1)^k \frac{\sigma_k}{k}(s-1)^k + \log\Gamma\left(1+\frac{s}{2}\right) + \log\pi^{s/2}\right]$$

As above we see that

$$Y_n(a_1+b_1+c_1,\ldots,a_n+b_n+c_n) = \sum_{k=0}^{n}\binom{n}{k}Y_{n-k}(a_1,\ldots,a_{n-k})\sum_{j=0}^{k}\binom{k}{j}Y_{k-j}(b_1,\ldots,b_{k-j})Y_j(c_1,\ldots,c_j)$$

but it is not clear whether the above methodology will result in an easily manipulated formula.

## 6. A formula for the Stieltjes constants in terms of the (exponential) complete Bell polynomials containing the eta constants $\eta_n$ as the arguments

We refer to (4.1)

$$\frac{\varsigma'(s)}{\varsigma(s)} = -\frac{1}{s-1} - \sum_{k=0}^{\infty} \eta_k (s-1)^k$$

and it is easily seen that this may equivalently be written as

$$\frac{d}{ds}[(s-1)\varsigma(s)] = -(s-1)\varsigma(s)\sum_{k=0}^{\infty} \eta_k (s-1)^k$$

As noted by Coffey [9] for differentiable functions $f(s)$ and $g(s)$ such that

$$f'(s) = g(s)f(s)$$

then we have in terms of the (exponential) complete Bell polynomials

$$f^{(n+1)}(s) = f(s)Y_{n+1}\left(g(s), g'(s), \ldots, g^{(n)}(s)\right)$$

In this particular case we have

$$f(s) = (s-1)\varsigma(s) \quad f(1) = 1 \quad \text{and} \quad g(s) = -\sum_{k=0}^{\infty} \eta_k (s-1)^k \quad g(1) = -\eta_0 = \gamma$$

$$g^{(j)}(s) = -\sum_{k=0}^{\infty} \eta_k k(k-1)\cdots(k-j+1)(s-1)^{k-j}$$

$$g^{(j)}(1) = -j!\eta_j$$

We recall from (2.10) that

$$\left.\frac{d^{n+1}}{ds^{n+1}}[(s-1)\varsigma(s)]\right|_{s=1} = (-1)^n (n+1)\gamma_n$$

and we then obtain the relationship for $n \geq 0$

(6.1)    $(-1)^n (n+1)\gamma_n = Y_{n+1}\left(\gamma, -1!\eta_1, \ldots, -n!\eta_n\right)$

For example with $n = 0$ and using (3.24) $Y_1(x_1) = x_1$ we recover $\gamma_0 = Y_1(\gamma) = \gamma$. Similarly with $n = 1$ we have using $Y_2(x_1, x_2) = x_1^2 + x_2$

$$-2\gamma_1 = Y_2(\gamma, -1!\eta_1) = \gamma^2 - \eta_1$$

which gives us

$$\eta_1 = \gamma^2 + 2\gamma_1$$

As a third example, using $Y_3(x_1, x_2, x_3) = x_1^3 + 3x_1 x_2 + x_3$ we obtain

$$3\gamma_2 = Y_3(\gamma, -1!\eta_1, -2!\eta_2) = \gamma^3 - 3\gamma\eta_1 - 2\eta_2$$

$$2\eta_2 = \gamma^3 - 3\gamma\eta_1 - 3\gamma_2 = \gamma^3 - 3\gamma(\gamma^2 + 2\gamma_1) - 3\gamma_2$$

and therefore we have

$$3\gamma_2 = -2\gamma^3 - 6\gamma\gamma_1 - 2\eta_2$$

We may write (6.1) as

$$(-1)^n (n+1)\gamma_n = Y_{n+1}\left(-0!\eta_0, -1!\eta_1, \ldots, -n!\eta_n\right)$$

Using the recurrence relation [28]

$$Y_{n+1}(x_1, \ldots, x_{n+1}) = \sum_{k=0}^{n} \binom{n}{k} Y_{n-k}(x_1, \ldots, x_{n-k}) x_{k+1} = \sum_{k=0}^{n} \binom{n}{k} Y_k(x_1, \ldots, x_k) x_{n-k+1}$$

$$= x_{n+1} + \sum_{k=1}^{n} \binom{n}{k} Y_k(x_1, \ldots, x_k) x_{n-k+1}$$

where $x_r = -(r-1)!\eta_{r-1}$ We see that

$$(-1)^n (n+1)\gamma_n = -n!\eta_n + \sum_{k=1}^{n} \binom{n}{k} (-1)^k k\gamma_{k-1}(n-k)!\eta_{n-k}$$

and hence we recover (4.4).

From (1.7) we have

$$\xi'(s) = -\xi(s) \sum_{k=1}^{\infty} (-1)^k \frac{\sigma_k}{k} k(s-1)^{k-1}$$

and in this case we designate $g(s) = -\sum_{k=1}^{\infty} (-1)^k \frac{\sigma_k}{k} k(s-1)^{k-1}$

$$g^{(j)}(s) = -\sum_{k=0}^{\infty} (-1)^k \frac{\sigma_k}{k} k(k-1)\cdots(k-j)(s-1)^{k-1-j}$$

$$g^{(j)}(1) = (-1)^j j! \sigma_{j+1}$$

$$\xi^{(n)}(s) = \xi(s) Y_n\left(g(s), g'(s), ..., g^{(n-1)}(s)\right)$$

$$\xi^{(n)}(1) = \frac{1}{2} Y_n\left(\sigma_1, -1!\sigma_2, ..., (-1)^{n-1}(n-1)!\sigma_n\right)$$

We have

$$\xi^{(n+1)}(1) = \frac{1}{2} Y_n\left(\sigma_1, -1!\sigma_2, ..., (-1)^n n!\sigma_{n+1}\right)$$

and using

$$Y_{n+1}(x_1, ..., x_{n+1}) = x_{n+1} + \sum_{k=1}^{n} \binom{n}{k} Y_k(x_1, ..., x_k) x_{n-k+1}$$

we obtain

(6.2) $$\xi^{(n+1)}(1) = \frac{1}{2}(-1)^n n! \sigma_{n+1} + \sum_{k=1}^{n} \binom{n}{k} (-1)^{n-k+1} \xi^{(k)}(1)(n-k+1)! \sigma_{n-k}$$

## Appendix A

### A brief survey of the (exponential) complete Bell polynomials

The (exponential) complete Bell polynomials may be defined by $Y_0 = 1$ and for $n \geq 1$

(A.1) $$Y_n(x_1, ..., x_n) = \sum_{\pi(n)} \frac{n!}{k_1! k_2! ... k_n!} \left(\frac{x_1}{1!}\right)^{k_1} \left(\frac{x_2}{2!}\right)^{k_2} ... \left(\frac{x_n}{n!}\right)^{k_n}$$

where the sum is taken over all partitions $\pi(n)$ of $n$, i.e. over all sets of integers $k_j$ such that

(A.2) $k_1 + 2k_2 + 3k_3 + ... + nk_n = n$

The complete Bell polynomials have integer coefficients and, by way of illustration, the first six are set out below [14, p.307]

(A.3) $$Y_1(x_1) = x_1$$

$$Y_2(x_1, x_2) = x_1^2 + x_2$$

$$Y_3(x_1, x_2, x_3) = x_1^3 + 3x_1 x_2 + x_3$$

$$Y_4(x_1, x_2, x_3, x_4) = x_1^4 + 6x_1^2 x_2 + 4x_1 x_3 + 3x_2^2 + x_4$$

$$Y_5(x_1, x_2, x_3, x_4, x_5) = x_1^5 + 10x_1^3 x_2 + 10x_1^2 x_3 + 15x_1 x_2^2 + 5x_1 x_4 + 10x_2 x_3 + x_5$$

$$Y_6(x_1, x_2, x_3, x_4, x_5, x_6) = x_1^6 + 6x_1 x_5 + 15x_2 x_4 + 10x_3^2 + 15x_1^2 x_4 + 15x_2^3 + 60x_1 x_2 x_3$$

$$+ 20x_1^3 x_3 + 45x_1^2 x_2^2 + 15x_1^4 x_1 + x_6$$

The complete Bell polynomials are also given by the exponential generating function in Comtet's book [14, p.134]

(A.4) $$\exp\left(\sum_{j=1}^{\infty} x_j \frac{t^j}{j!}\right) = 1 + \sum_{n=1}^{\infty} Y_n(x_1,\ldots,x_n) \frac{t^n}{n!} = \sum_{n=0}^{\infty} Y_n(x_1,\ldots,x_n) \frac{t^n}{n!}$$

Let us now consider a function $f(x)$ which has a Taylor series expansion around $x$: we have

$$e^{f(x+t)} = \exp\left(\sum_{j=0}^{\infty} f^{(j)}(x) \frac{t^j}{j!}\right) = e^{f(x)} \exp\left(\sum_{j=1}^{\infty} f^{(j)}(x) \frac{t^j}{j!}\right)$$

$$= e^{f(x)} \left\{ 1 + \sum_{n=1}^{\infty} Y_n\left(f^{(1)}(x), f^{(2)}(x), \ldots, f^{(n)}(x)\right) \frac{t^n}{n!} \right\}$$

We see that

$$\frac{d^m}{dx^m} e^{f(x)} = \frac{\partial^m}{\partial x^m} e^{f(x+t)} \bigg|_{t=0} = \frac{\partial^m}{\partial t^m} e^{f(x+t)} \bigg|_{t=0}$$

and we therefore obtain (as noted by Kölbig [22] and Coffey [9])

(A.5) $$\frac{d^m}{dx^m} e^{f(x)} = e^{f(x)} Y_m\left(f^{(1)}(x), f^{(2)}(x), \ldots, f^{(m)}(x)\right)$$

Differentiating (A.5) we see that

$$\frac{d}{dx} Y_m\left(f^{(1)}(x), f^{(2)}(x),..., f^{(m)}(x)\right) = Y_{m+1}\left(f^{(1)}(x), f^{(2)}(x),..., f^{(m+1)}(x)\right)$$

$$- f^{(1)}(x) Y_m\left(f^{(1)}(x), f^{(2)}(x),..., f^{(m)}(x)\right)$$

As an example of (A.5), letting $f(x) = \log \Gamma(x)$ we obtain

(A.6) $$\frac{d^m}{dx^m} e^{\log \Gamma(x)} = \Gamma^{(m)}(x) = \Gamma(x) Y_m\left(\psi(x), \psi^{(1)}(x),..., \psi^{(m-1)}(x)\right)$$

$$= \int_0^\infty t^{x-1} e^{-t} \log^m t \, dt$$

and since [30, p.22]

$$\psi^{(p)}(x) = (-1)^{p+1} p! \varsigma(p+1, x)$$

we may express $\Gamma^{(m)}(x)$ in terms of $\psi(x)$ and the Hurwitz zeta functions. In particular, Coffey [9] notes that

(A.7) $$\Gamma^{(m)}(1) = Y_m(-\gamma, x_1,..., x_{m-1})$$

where $x_p = (-1)^{p+1} p! \varsigma(p+1)$.

Values of $\Gamma^{(m)}(1)$ are reported in [30, p.265] for $m \leq 10$ and the first three are

$$\Gamma^{(1)}(1) = -\gamma$$

$$\Gamma^{(2)}(1) = \varsigma(2) + \gamma^2$$

$$\Gamma^{(3)}(1) = -[2\varsigma(3) + 3\gamma\varsigma(2) + \gamma^3]$$

The general form is

$$\Gamma^{(m)}(1) = (-1)^m \sum_{j=1}^m \varepsilon_{mj}$$

where $\varepsilon_{mj}$ are positive constants.

$\square$

The following is extracted from an interesting series of papers written by Snowden [29, p.68].

Let us consider the function $f(x)$ with the following Maclaurin expansion

(A.8) $$\log f(x) = b_0 + \sum_{n=1}^{\infty} \frac{b_n}{n} x^n$$

and we wish to determine the coefficients $a_n$ such that

(A.9) $$f(x) = \sum_{n=0}^{\infty} a_n x^n$$

By differentiating (A.8) and multiplying the two power series, we get

(A.10) $$na_n = \sum_{k=1}^{n} b_k a_{n-k}$$

Upon examination of this recurrence relation Snowden reports it is easy to see that

(A.11) $$n! a_n = a_0 [b_1, -b_2, b_3, ..., (-1)^{n+1} b_n]$$

where the symbol $[c_1, c_2, c_3, ..., c_n]$ is defined as the $n \times n$ determinant

$$\begin{vmatrix} c_1 & c_2 & c_3 & c_4 & . & . & . & c_n \\ (n-1) & c_1 & c_2 & c_3 & . & . & . & c_{n-1} \\ 0 & (n-2) & c_1 & c_2 & . & . & . & c_{n-2} \\ 0 & 0 & (n-3) & c_1 & . & . & . & c_{n-3} \\ \vdots & \vdots & \vdots & \vdots & \vdots & \vdots & \vdots & \vdots \\ 0 & 0 & 0 & 0 & 0 & 0 & 1 & c_1 \end{vmatrix}$$

Since $\log f(0) = \log a_0 = b_0$ we have

(A.12) $$f(x) = e^{b_0} \left[ 1 + \sum_{n=1}^{\infty} [b_1, -b_2, b_3, ..., (-1)^{n+1} b_n] \frac{x^n}{n!} \right]$$

Multiplying (A.9) by $\alpha$ it is easily seen that (after correcting the misprint in Snowden's paper)

(A.13) $$f^{\alpha}(x) = e^{\alpha b_0}\left[1 + \sum_{n=1}^{\infty}[\alpha b_1, -\alpha b_2, \alpha b_3, ..., (-1)^{n+1}\alpha b_n]\frac{x^n}{n!}\right]$$

and, in particular, with $\alpha = -1$ we obtain

(A.14) $$\frac{1}{f(x)} = e^{-b_0}\left[1 + \sum_{n=1}^{\infty}[-b_1, b_2, -b_3, ..., (-1)^n b_n]\frac{x^n}{n!}\right]$$

Differentiating (A.13) with respect to $\alpha$ would give us an expression for $f^{\alpha}(x)\log f(x)$.

Differentiating (A.8) gives us

$$f'(x) = f(x)\sum_{k=1}^{\infty} b_k x^{k-1} = f(x)g(x)$$

and hence

$$f^{(n)}(x) = f(x)Y_n\left(g(x), ..., g^{(n-1)}(x)\right)$$

Alternatively, differentiating (A.9) results in

$$f^{(n)}(x) = \sum_{k=0}^{\infty} a_k k(k-1)\cdots(k-n+1)x^{k-n}$$

and letting $x = 0$ gives us

(A.15) $$n!a_n = e^{b_0} Y_n\left(b_1, 1!b_2, ..., (n-1)!b_n\right)$$

From (A.11) we then see that

(A.16) $$Y_n\left(b_1, 1!b_2, ..., (n-1)!b_n\right) = [b_1, -b_2, b_3, ..., (-1)^{n+1} b_n]$$

and more generally

(A.17) $$Y_n(x_1, ..., x_n) = \left[\frac{x_1}{0!}, -\frac{x_2}{1!}, ..., (-1)^{n+1}\frac{x_n}{(n-1)!}\right]$$

We now write (A.13) as

$$f^{\alpha}(x) = e^{\alpha b_0}\left[1 + \sum_{n=1}^{\infty} D_n(\alpha)\frac{x^n}{n!}\right]$$

where for convenience we have designated $D_n(\alpha)$ by

$$D_n(\alpha) = [\alpha b_1, -\alpha b_2, \alpha b_3, ..., (-1)^{n+1} \alpha b_n]$$

Letting $\alpha = 0$ we see that

$$1 = 1 + \sum_{n=1}^{\infty} D_n(0) \frac{x^n}{n!}$$

and therefore we have

$$D_n(0) = 0$$

which is only to be expected since all of the elements of the top row of the defining determinant are equal to zero.

Differentiation with respect to $\alpha$ results in

$$f^{\alpha}(x) \log f(x) = e^{\alpha b_0} \sum_{n=1}^{\infty} D'_n(\alpha) \frac{x^n}{n!} + b_0 e^{\alpha b_0} \left[ 1 + \sum_{n=1}^{\infty} D_n(\alpha) \frac{x^n}{n!} \right]$$

Letting $\alpha = 0$ gives us

$$\log f(x) = \sum_{n=1}^{\infty} D'_n(0) \frac{x^n}{n!} + b_0 \left[ 1 + \sum_{n=1}^{\infty} D_n(0) \frac{x^n}{n!} \right]$$

$$= b_0 + \sum_{n=1}^{\infty} D'_n(0) \frac{x^n}{n!}$$

We started out with (A.8)

$$\log f(x) = b_0 + \sum_{n=1}^{\infty} \frac{b_n}{n} x^n$$

and, hence equating coefficients in the two power series, we obtain

$$D'_n(0) = (n-1)! b_n$$

We have

$$\frac{d^m}{dx^m} f^\alpha(x) = e^{\alpha b_0} \sum_{n=1}^{\infty} D_n(\alpha) \frac{n(n-1)\cdots(n-m+1)x^{n-m}}{n!}$$

$$\frac{d}{d\alpha}\frac{d^m}{dx^m} f^\alpha(x) = e^{\alpha b_0} \sum_{n=1}^{\infty} D'_n(\alpha) \frac{n(n-1)\cdots(n-m+1)x^{n-m}}{n!}$$

$$+ b_0 e^{\alpha b_0} \sum_{n=1}^{\infty} D_n(\alpha) \frac{n(n-1)\cdots(n-m+1)x^{n-m}}{n!}$$

and hence

$$\left.\frac{d}{d\alpha}\frac{d^m}{dx^m} f^\alpha(x)\right|_{x=1,\alpha=0} = \sum_{n=1}^{\infty} D'_n(0) \frac{n(n-1)\cdots(n-m+1)}{n!}$$

We then have

$$\left.\frac{d}{d\alpha}\frac{d^m}{dx^m} f^\alpha(x)\right|_{x=1,\alpha=0} = \sum_{n=1}^{\infty} \frac{b_n}{n} n(n-1)\cdots(n-m+1)$$

which is equivalent to

$$\left.\frac{d^m}{ds^m} \log f(x)\right|_{x=1} = \sum_{n=1}^{\infty} \frac{b_n}{n} n(n-1)\cdots(n-m+1)$$

Donal F. Connon
Elmhurst
Dundle Road
Matfield
Kent TN12 7HD
dconnon@btopenworld.com